\documentclass[a4paper,12pt]{article}

\usepackage{citesort,geometry,enumerate,amsmath,amssymb,latexsym,theorem}
\usepackage{epic,eepic} 
\usepackage{xypic}
\xyoption{curve}

{\theorembodyfont{\rmfamily}\newtheorem{example}{Example}[section]}
{\theorembodyfont{\rmfamily}\newtheorem{defn}[example]{Definition}}
\newtheorem{prop}[example]{Proposition}
\newtheorem{thm}[example]{Theorem}
{\theorembodyfont{\rmfamily}}
\newtheorem{cor}[example]{Corollary}
\newtheorem{lem}[example]{Lemma}

\def\<{\langle}
\def\>{\rangle}

\newcommand{\QED}{\hfill \rule[.25mm]{1.9mm}{0mm}}
\newenvironment{pf}{{\bf Proof:} }{\hfill $\Box$

\mbox{}}

\def\geq{\geqslant}
\def\Aut{\mathsf{Aut}}
\def\C{\mathcal{C}}
\def\D{\mathcal{D}}
\def\eps{\varepsilon}

\def\Do{\mathsf{D}}
\def\epsilon{\varepsilon}

\newcommand{\directs}[2]{\def\objectstyle{\scriptstyle}  \objectmargin={0pt}
\xy
(0,4)*+{}="a",(0,-2)*+{\rule{0em}{1.5ex}#2}="b",(7,4)*+{\;#1}="c"
\ar@{->} "a";"b" \ar @{->}"a";"c" \endxy }
\newcommand{\xdirects}[2]{\def\objectstyle{\scriptstyle}  \objectmargin={0pt}
\xy
(0,0)*+{}="a",(0,-6)*+{\rule{0em}{1.5ex}#2}="b",(7,0)*+{\;#1}="c"
\ar@{->} "a";"b" \ar @{->}"a";"c" \endxy }

\begin{document}

\title{\huge{\bf Towards a 2-dimensional notion of holonomy}
\thanks{MSC-class:
18D05, 58H05, 22A22, 22E99. Key words: crossed module, homotopy,
 double groupoid, linear section,
 Lie groupoids, holonomy.}}

\small{ \author  {Ronald Brown  \\ University of Wales  \\ School
of Mathematics \\ Bangor, Gwynedd  LL57 1UT,\\ U.K.\\E-mail:
r.brown@bangor.ac.uk   \\ \and \.{I}lhan \.{I}\c{c}en    \\
University of  \.{I}n\"{o}n\"{u} \\ Faculty of Science and Art
\\ Department of
Mathematics
\\ Malatya/ Turkey\\E-mail: iicen@inonu.edu.tr}
\maketitle

\begin{center}
University of Wales Bangor Maths Preprint 00.14
\end{center}

\begin{abstract}
Previous work (Pradines, 1966, Aof and Brown, 1992) has given a
setting for a holonomy Lie groupoid of a locally Lie groupoid.
Here we develop analogous 2-dimensional notions starting from a
locally Lie crossed module of groupoids. This involves replacing
the Ehresmann notion of a local smooth coadmissible section of a
groupoid by a local smooth coadmissible homotopy (or free
derivation) for the crossed module case. The development also has
to use corresponding notions for certain types of double
groupoids. This leads to  a holonomy Lie groupoid rather than
double groupoid, but one which involves the $2$-dimensional
information.
\end{abstract}
\addcontentsline{toc}{chapter}{\protect\numberline{
}{Introduction}}

\newpage

\section*{Introduction}
An intuitive notion of holonomy is that it is concerned with
\begin{quote}
{\em Iterations of local procedures}
\end{quote}
and it is detected when such an iteration returns to its starting
point but with a `change of phase'. An aim of theoretical
developments is to provide  appropriate algebraic settings for
this intuition.  Major areas in which holonomy occurs are in
foliation theory and in differential geometry, and it is with
ideas arising out of the first area with which we are mainly
concerned.

The leaves of a  foliation on a manifold $M$ are path connected
subsets of $M$ forming a partition of $M$, and so an equivalence
relation $R$, say, on $M$. This equivalence relation is usually
not a submanifold of $M$ -- for example in the classical foliation
of the M\"obius Band, $R$ is of dimension 3 but has self
intersections. However the local structure describing the
foliation determines a subset $W$ of $R$ such that $W$ is a smooth
manifold. Further, $R$ can be regarded as a groupoid with the
usual multiplication $(x,y)(y,z)=(x,z)$ for $(x,y),(y,z)\in R$.
The pair $(R,W)$ becomes what is known as a {\em locally Lie
groupoid} -- that is, the groupoid operations are as smooth as
they can be on $W$, granted  that they are not totally defined on
$W$. This concept was first formulated by J. Pradines in 1966
\cite{Pr1} and called `un morceau d'un groupo\"\i de
diff\'erentiables'. This description of $(R,W)$ for a foliation
was due to Kubarski in 1987 \cite{Ku} and independently to Brown
and Mucuk in 1995 \cite{Br-Muc1}, though with slightly differing
conditions.

It is classical that if $G$ is a group and $W$ is a subset of $G$
containing $1$ and $W$ has a topology, then under reasonable
conditions on the pair $(G,W)$ the topology of $W$ can be
transported around $G$ to give a base for a topology so that  $G$
becomes a topological group in which $W$ is an open subspace. We
say that a {\it locally Lie group is extendible}. This is not
generally true for locally Lie groupoids -- and the above $(R,W)$
provides a counterexample. A basic reason for this is that while
for a topological group $G$ left multiplication $L_g$ by an
element $g$ maps open sets to open sets, this is no longer true
for general topological groupoids, since the domain of $L_g$ is
usually not open in $G$.

Ehresmann realised that this failure of this useful property of
left multiplication by an element  could be remedied by replacing
the element $g$ by a local smooth coadmissible section through $g$
(for technical reasons we are replacing `admissible' by
`coadmissible'). This is a continuous function $s:U\to G$ for some
open subset of the base space $X$ and such that $\beta s=1_U$, and
$\alpha s$ maps $U$ homeomorphically to an open subset of $X$
(here $\alpha, \beta$ are source and target maps of $G$). Then
$L_s$ is the partial map on $G$ given by $h\mapsto s(\alpha h)+
h$, and $L_s$ does indeed map open sets of $G$ to open sets.
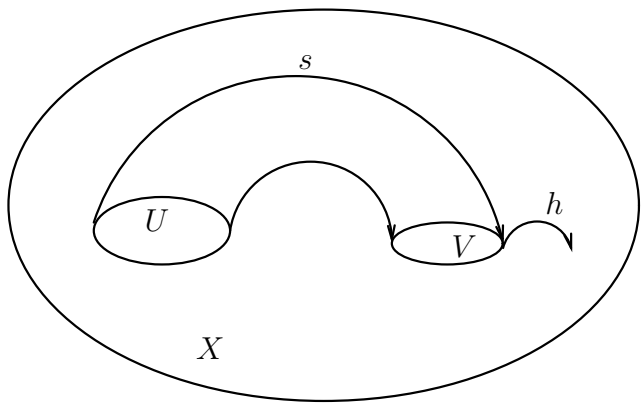
\begin{center}

\setlength{\unitlength}{0.0006in}
\begingroup\makeatletter\ifx\SetFigFont\undefined
\def\x#1#2#3#4#5#6#7\relax{\def\x{#1#2#3#4#5#6}}%
\expandafter\x\fmtname xxxxxx\relax \def\y{splain}%
\ifx\x\y   
\gdef\SetFigFont#1#2#3{%
  \ifnum #1<17\tiny\else \ifnum #1<20\small\else
  \ifnum #1<24\normalsize\else \ifnum #1<29\large\else
  \ifnum #1<34\Large\else \ifnum #1<41\LARGE\else
     \huge\fi\fi\fi\fi\fi\fi
  \csname #3\endcsname}%
\else \gdef\SetFigFont#1#2#3{\begingroup
  \count@#1\relax \ifnum 25<\count@\count@25\fi
  \def\x{\endgroup\@setsize\SetFigFont{#2pt}}%
  \expandafter\x
    \csname \romannumeral\the\count@ pt\expandafter\endcsname
    \csname @\romannumeral\the\count@ pt\endcsname
  \csname #3\endcsname}%
\fi \fi\endgroup
\begin{picture}(5566,3481)(0,-10)
\thicklines
\path(3347.385,1551.455)(3383.000,1433.000)(3407.319,1554.279)
\put(2666.750,1399.250){\arc{1434.089}{3.2938}{6.2361}}
\path(4899.791,1467.141)(4958.000,1358.000)(4958.000,1481.693)
\put(4658.000,1283.000){\arc{618.466}{3.3866}{6.0382}}
\put(2783,1733){\ellipse{5550}{3450}}
\put(1358,1508){\ellipse{1200}{600}}
\path(4301.236,1542.899)(4358.000,1433.000)(4359.631,1556.682)
\put(2536.966,1003.172){\arc{3742.147}{3.4567}{6.0514}}
\put(3870,1396){\ellipse{974}{374}}
\put(4733,1658){\makebox(0,0)[lb]{\smash{{{\SetFigFont{12}{14.4}{rm}$h$}}}}}
\put(1208,1508){\makebox(0,0)[lb]{\smash{{{\SetFigFont{12}{14.4}{rm}$U$}}}}}
\put(3908,1280){\makebox(0,0)[lb]{\smash{{{\SetFigFont{12}{14.4}{rm}$V$}}}}}
\put(1658,383){\makebox(0,0)[lb]{\smash{{{\SetFigFont{12}{14.4}{rm}$X$}}}}}
\put(2558,2933){\makebox(0,0)[lb]{\smash{{{\SetFigFont{12}{14.4}{rm}$s$}}}}}
\end{picture}
\end{center}

The local section $s$ of these is {\em through} an element $g\in
G$ if $s(\beta g)=g$. We adopt the idea  that $s$ is a kind of
thickening or localisation of $g$, that is it extends $g$ to a
local neighbourhood.

Thus  one part of the legacy of Charles Ehresmann \cite{Eh1} is
the realisation that from the point of view of the topology of a
Lie groupoid $G$ the interest is less in the elements but more in
the local coadmissible sections and their actions on $G$. In the
case of the locally Lie groupoid $(G,W)$, the local $\C^r$
coadmissible sections with values in $W$ can be regarded as `local
procedures'. They have a multiplication due to Ehresmann, and
germs of their iterates  form a groupoid $J^{r}(G,W)$. Pradines
realised that a quotient groupoid of $J^{r}(G,W)$ does obtain a
Lie groupoid structure with $W$ as an open subspace, and is the
minimal such overgroupoid of $G$. This is reasonably called the
{\em holonomy groupoid} of $(G,W)$ and written $Hol(G,W)$. It
encapsulates many of the  usual intuitions behind the holonomy
groupoid of a foliation (as described for example in Winkelnkemper
\cite{Wi} and Phillips \cite{Ph}) but in much greater generality.
The major lines of Pradines' construction were explained to the
first author in the period 1981-87 and were published in full
detail by Aof and Brown in \cite{Aof-Br} with the generous
agreement of Pradines.

It is this intuition  which we would like to extend  to dimension
$2$, using one of several notions of $2$-dimensional groupoid.

A key idea as to how this might work is that a coadmissible
section $s:X\to G$ of a groupoid $G$ on $X$ can also be regarded
as a homotopy $s:\Delta(s)\simeq1$ from an automorphism
$\Delta(s)$ of $G$ to the identity on $G$. Thus the construction
of $Hol(G, W)$ can be expected to be realised in situations where
we have a notion of homotopy. In dimension $2$, this is available
in the most well worked out way  for {\em crossed modules of
groupoids} and this is one of the major areas in which we work.
These objects are in fact equivalent to {\em edge symmetric double
groupoids with connection}, and also to {\em $2$-groupoids}. We
find it possible to work with the first two of these structures
and these give the framework for our constructions.

We work in a kind of inductive and non symmetric format: that is
we suppose given a Lie structure in one of the two dimensions and
use a localisation at the next dimension -- precise definitions
are given later.

Applications are expected to be in  situations with multiple
geometric structures, such as foliated bundles. Experience in this
area of `higher dimensional group theory' has suggested that it is
necessary to build first a strong feel for the appropriate
algebra. For example, the notion of double groupoid was considered
by the first author in 1967; double groupoids  with connection
were found by Brown and Spencer in 1971, but the homotopical
notion of the homotopy double groupoid of a pair of spaces, and
its application to a $2$-dimensional Van Kampen Theorem for the
crossed module $\pi_2(X, A)\to \pi_1(A)$ was not found by Brown
and Higgins till 1974, and published only in 1978. Once the
algebra was developed and linked with the geometry, then quite
novel geometrical  results were obtained.

By following this paradigm, we intend to come nearer to
2-dimensional extensions of the notions of transport along a path.
This should give ideas of, for example, transport over a surface,
and pave the way for further extensions to all dimensions. It is
hoped that this will lead to a deeper understanding of higher
dimensional constructions and operations in differential topology.
It is for these reasons, that we still title this paper with the
word `towards'.

This work is also intended to be a continuation of other work
applying double groupoids in differential topology such as
Mackenzie~\cite{Mac2}, Brown-Mackenzie~\cite{Br-Mac}.

In the first section, we outline our plan of the construction.

\section{Plan of the work}

For a 2-dimensional version, there are a number of possible
choices for analogues of groupoids, for example double groupoids,
$2$-groupoids, crossed modules of groupoids. We are not able at
this stage to give a version of holonomy for the most general
locally Lie double groupoids. It seems reasonable therefore to
restrict attention to those forms of double groupoids whose
algebra is better  understood, and we therefore considered the
possibility of a theory for one of the equivalent categories
\[ (CrsMod) \sim (2-Grpd)\sim  (DGrpd!), \]
which denote respectively the categories of crossed modules over
groupoids, $2$-groupoids and ``edge symmetric  double groupoids
with connection''. For consideration of homotopies the outside
seem more convenient. There are more general versions of double
groupoids (see Brown and Mackenzie \cite{Br-Mac}) whose
consideration we leave to further work.

The steps that we take  are as follows:

(i) We need to formulate the notion of a locally Lie structure on
a double groupoid ${\cal D}({\cal C})$ which corresponds to a
crossed module ${\cal C} = ( C, G, \delta)$ with base space $X$.
For this reason, here $(G, X)$ is supposed a Lie groupoid and a
smooth manifold structure on a set $W$ such that $X\subseteq
W\subseteq C.$ Then  $({\cal D}({\cal C}), W^G)$ can given as a
locally Lie groupoid over $G$, where  $W^G$ is the subset of
${\cal D}({\cal C})$ given by   $$\{
  \left( \spreaddiagramrows{-2pc}  \spreaddiagramcolumns{-1.5pc}
\objectmargin{0.05pc} w_1 : \diagram  & d & \\
                                     b & & c \\ & a &   \enddiagram \right)
:  w_1\in W, a, b, c\in G, \beta (b) = \alpha (a), \beta (a) =
\beta (c) = \beta (w_1), \ d=b+a+\delta(w_1)-c  \} .$$

(ii) Next, we are replacing local coadmissible sections of a
groupoid by local linear coadmissible sections of an edge
symmetric double groupoid. We define a product on the set of all
local linear coadmissible sections. This easily leads to a
2-dimensional version of  $\Gamma({\cal D}({\cal C}))$ of $\Gamma
(G)$, again  an inverse semigroup.

(iii) Now we form germs  $[s]_a$ of $s$, where $a\in G$, $s\in
\Gamma ({\cal D}({\cal C}))$. The set of these germs forms a
groupoid $J({\cal D}({\cal C}))$ over $G$.

(iv)  A key matter for decision is that of the final map $\psi $
and its values on $[s]_a$, since the formation of the holonomy
groupoid will involve $Ker\; \psi$.

(v) A further question is that of deciding the meaning of the
generalisation   to dimension 2 of the term ``enough local linear
coadmissible sections''. This requires further discussion.

Recall from \cite{Aof-Br} that, in the groupoid case, we  ask that
for any $a\in G$ there is a local coadmissible section $s$ such
that $\beta a\in \Do(s)$ (where $\Do(s)$ is the domain of $s$) and
$s \beta a = a.$ Under certain conditions, we require $s$ to be
smooth and such that $\alpha s$ is a diffeomorphism of open sets.
The intuition here is that $a\in G$ can be regarded as a
deformation of $\beta a$, and $s$ gives a ``thickening'' of this
deformation.

In dimension 2, we therefore suppose given $a\in G(x, y)$ and
$b\in G(z, x)$, $c\in G(w, y)$ and $w_1\in C(y)$. $$
\spreaddiagramrows{2pc} \spreaddiagramcolumns{2pc}
\def\objectstyle{\ssize} \def\labelstyle{\ssize}
\def\labelstyle{\textstyle}
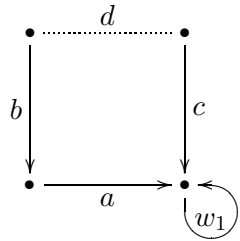
\diagram
\bullet \dto_{b}         \rdotted^{d}  & \bullet \dto^{c}   \\
\bullet^{}     \rto_{a}   & \bullet  \todr^{w_1}
\enddiagram
$$ where $d=b+a+\delta (w_1)-c$.  Then a local linear coadmissible
section $s=(s_0,s_1)$ will be ``through $ w = \left(
\spreaddiagramrows{-2pc} \spreaddiagramcolumns{-1.5pc}
\objectmargin{0.05pc} w_1 : \diagram & d & \\
                                     b & & c \\ & a &   \enddiagram \right)$
''  if $s_0 x = b, s_0 y = c$ and $s_1 a = w_1$. Our ``final map''
$\psi$ will be a morphism from $J ({\cal D}({\cal C}))$  to a
groupoid. This groupoid ${\cal D}({\cal C})$ will be one of the
groupoid structures of the double groupoid associated to the
crossed module ${\cal C} = {(C, G, \delta )}$. We write $$\psi
([s]_a) = s(a)  = \left( \spreaddiagramrows{-2pc}
\spreaddiagramcolumns{-1.5pc}
   \objectmargin{0.05pc} s_1(a) : \diagram  & f_1(a) & \\
                                   s_0(x) &         & s_0(y) \\ & a &   \enddiagram \right),$$
so that  the value of $\psi $ on $[s]_a$ does use all the
information given by $s = ( s_0, s_1)$ at the arrow $a\in G$. This
explains why our theory  develops crossed modules and double
groupoids in parallel.


\section{Crossed modules and  edge symmetric double groupoids\\ with connection}

In previous paper \cite{Br-Ic2}, we  explored  the idea that a
natural generalisation to crossed modules of the notion of
coadmissible section for groupoids is that of {\bf coadmissible
homotopy}. This arises naturally from the work of Brown-Higgins
\cite{Br-Hig2} on homotopies for crossed complexes over groupoids.

We recall the definition of crossed modules of groupoids. The
original reference is  Brown-Higgins \cite{Br-Hig4}, but see also
\cite{Br-Ic2}.

The source and target maps of a groupoid $G$ are written $\alpha,
\beta$ respectively. If  $G$ is {\em totally intransitive}, i.e.
if $\alpha=\beta$, then we usually use the notation $\beta$. The
composition in a groupoid $G$ of elements $a,b$ with $\beta a =
\alpha  b$ will be written additively, as $a + b$.  The main
reason for this is the convenience for dealing with combinations
of inverses and actions.

\begin{defn}
Let $G, C$ be groupoids over the same object set and let $C$ be
totally intransitive. Then an {\bf action}  of $G$ on $C$ is given
by a partially defined function
\[ \xymatrix {       C\times G \ar@{o->}[r]& C }        \]
written $(c, a)\mapsto c^a$, which satisfies

\begin{enumerate}[(i)]
  \item $c^a$ is defined if and only if $\beta (c) = \alpha (a) $, and then
$\beta (c^a) = \beta (a)$,
  \item  $(c_1 + c_2)^a = {c_1}^a + {c_2}^a$,
  \item $c_1^{a+b} = (c_1^{a})^{b}$ and  $c_1 ^{e_x} = c_1$
\end{enumerate}
for all $c_1, c_2 \in C(x, x)$, $a\in G(x, y) $, $b\in G(y, z)$.
\end{defn}

\begin{defn}
A {\bf crossed  module of groupoids} \cite{Br-Hig4} consists of a
morphism $\delta : C\rightarrow G$  of groupoids $C$ and $G$ which
is the identity on the object sets such that $C$ is totally
intransitive, together with an action of $G$ on $C$ which
satisfies
\begin{enumerate}[(i)]
  \item $\delta (c^a ) = -a + \delta c +a$,
  \item  $c^{\delta c_1} = -{c_1}+ c + c_1 $,
\end{enumerate}
for $c, c_1\in C(x, x)$, $a\in G(x, y)$.
\end{defn}

\begin{defn} \cite{Br-Ic2}
Let ${\cal C} =(C, G, \delta )$ a crossed module with the base
space $X$. A {\bf free derivation} $s$ is a pair of   maps $s_0
\colon X\to G,$ $s_1:G\to C$ which satisfy the following
\begin{align*}
\beta (s_0 x)          &=  x,  \ \ x\in X \\ \beta (s_1 a) &=
\beta (a),\ \ a\in G,         \\ s_1(a + b ) &=  s_1( a)^{b} +
s_1(b), \ \  a, b\in G.
\end{align*}
\end{defn}
Let ${\sf {\sf FDer}}({\cal C})$ be the set of free derivations of
${\cal C}$.
\\ We proved in a previous paper \cite{Br-Ic2} that if $s$ is a free
derivation of the crossed module ${\cal C} =(C, G, \delta )$ over
groupoids, then the formulae
\begin{align*}
f_0(x)  &=  \alpha s_0 (x), x \in X\\ f_1(a)  & = s_0(\alpha a) +
a +\delta s_1(a) -s_0(\beta a), a \in G\\ f_2(c) &=  (c+s_1\delta
c)^{-s_0\beta (c)}, c \in C
\end{align*}
define an endomorphism  $\Delta(s)= (f_0, f_1, f_2)$ of ${\cal
C}$.

Also  ${\sf FDer}({\cal C})$  has a monoid structure with the
following multiplication \cite{Br-Ic2}:
\begin{align*}
(s*t)_{\epsilon}(z) &= \begin{cases} t_1(z)+
(s_1g_1(z))^{t_0(\beta z)}, & \epsilon=1, z\in G(x, y), \\
(s_0g_0(z))+t_0(z), & \epsilon=0, \ \ z\in X, \ \
\end{cases}
\end{align*}
where $g = (g_0, g_1, g_2) = \Delta(t)$. This multiplication for
$\epsilon=0$  give us  Ehresmann's multiplication of coadmissible
sections \cite{Eh1}, and for  $\epsilon=1$ gives the
generalisation to free derivations of the  multiplication of
derivations introduced by Whitehead \cite{Wh1}.

Let $ {\sf FDer}^*({\cal C})$    denote the group of invertible
elements of this monoid. Then each element of ${\sf FDer}^*({\cal
C})$  is also called
 a {\bf coadmissible homotopy}.

\begin{thm}\label{14}
Let $s\in {\sf FDer}({\cal C})$  and let $f=\Delta (s)$. Then the
following conditions are equivalent: {\em
\begin{enumerate}[(i)]
  \item  $s\in {\sf FDer}^*({\cal C})$,
  \item   $f_1\in \Aut(G)$,
  \item   $f_2\in \Aut(C)$.\hfill$\Box$
\end{enumerate}}
\end{thm}
The proof is given in \cite{Br-Ic2}.

We  shall deal with also double groupoids especially edge
symmetric double groupoids. A double groupoid $\D$ is a groupoid
object internal to the category of groupoids. It may also be
represented as consisting of four groupoid structures
\[   (D, H, \alpha_1, \beta_1, + _1, \eps_1)
\ \ \ (D,V,\alpha_2, \beta_2, + _2, \eps_2) \] \[ (V, X, \alpha,
\beta,+,\eps)   \  \   \  (H, X, \alpha, \beta, + ,\eps)
\] as partially shown in the diagram
$$\diagram
D\rto<.5ex>\rto<-.5ex>\dto<.5ex>\dto<-.5ex>&V\dto<.5ex>\dto<-.5ex>\\
H\rto<.5ex>\rto<-.5ex>&X \enddiagram. $$ Here $H,\; V$ are called
the {\em horizontal} and {\em vertical edge groupoids}. The
functions written with $\alpha, \beta$ are the source and target
maps of the groupoids, and the $\eps$ denote the functions giving
the identity (zero) elements. Thus  ${\cal D}$ has two groupoid
structures $+_1,+_2$ over groupoids $H$ and $V$, which are
themselves groupoids on the common set $X$. These are all subject
to the compatibility condition that the structure maps of each
structure on ${\cal D}$ are morphisms with respect to the other.

Elements of $ D$ are pictured as squares $$
\spreaddiagramrows{2pc}
\spreaddiagramcolumns{2pc}
\def \objectstyle {\ssize} \def \labelstyle{\sssize}
\def\labelstyle{\textstyle}\xymatrix{
  \bullet \ar [d]_{\alpha_2(w)=v_1} \ar @{} [dr]|w      \ar
  [r]^{\alpha_1(w)=h_1}                &\bullet \ar [d]^{\beta_2(w)=v_2}  \\
  \bullet       \ar [r]_{\beta_1(w)=h_2}            &\bullet}
 \qquad \xdirects{2}{1} $$ in which $v_1, v_2\in V$ are the source
and target of $w$ with respect to the horizontal structure $+_2$
on $ D$, and $h_1, h_2\in H$ are the  source and target with
respect to the vertical structure $+_1$. (Thus our convention for
directions is the same as that for matrices.) Note that we use
addition for the groupoid compositions, but principally because it
makes for the easier use of negatives rather than inverses.
However,somewhat inconsistently, we shall tend use $1$ rather than
$0$ for the unit for $+$.  For further information on double
groupoids, we refer to \cite{Br-Sp1,Br3,Br-Mac,Br-Mo}.

Double groupoids  were introduced by Ehresmann in the early 1960's
\cite{Eh2,Eh3}, but primarily as instances of double categories,
and as a part of a general exploration of categories with
structure. Since that time their main use has been in homotopy
theory. Brown-Higgins \cite{Br-Hig3} gave the earliest example of
a ``higher homotopy groupoid'' by associating to a pointed pair of
spaces $(X, A)$ a edge symmetric double groupoid  with  connection
$\rho (X, A)$ in the sense of Brown and Spencer (see below). In
such a double groupoid, the vertical and horizontal edge
structures $H$ and $V$ coincide.  In terms of this functor $\rho$,
\cite{Br-Hig3} proved a Generalised Van Kampen Theorem, and
deduced from it a Van Kampen Theorem for the second relative
homotopy  crossed module  $\pi_2(X, A) \to \pi_1(A)$.

We shall initially be interested in edge symmetric double
groupoid, i.e. those in which the horizontal and vertical edge
groupoids coincide. These were called special double groupoids in
\cite{Br-Sp1}. In this case we write $G$ for $H=V$. The  main
result of Brown-Spencer in \cite{Br-Sp1} is that an edge symmetric
double groupoid with  connection (see later) whose double base is
a singleton is entirely determined by a certain crossed module it
contains; as explained above, crossed modules had arisen much
earlier in the work of J.H.C. Whitehead \cite{Wh2} on
2-dimensional homotopy. This result of Brown and Spencer  is
easily extended to give an equivalence between arbitrary edge
symmetric  double groupoids with connection and crossed modules
over groupoids; this is included in the result of \cite{Br-Hig4}.

We give this extended result as in \cite{Br-Hig4} and
\cite{Br-Mac}. The latter paper describes a more general class of
double groupoids and the two crossed modules given below form part
of the `core diagram' of a double groupoid described in
\cite{Br-Mac}.

The method which is used here can be found in \cite{Br-Sp1}. The
sketch proof is as follows:

Let ${\cal D} =(D, H, V, X)$ be a double groupoid. We show that
${\cal D}$ determines two crossed modules over groupoids.

Let $x\in X$ and let
\[  H(x) = \{ a\in H : \alpha (a) = \beta(a)=x \}. \]
We define $V(x)$ similarly. We put
\[ \Pi(D, H, x) = \{ w\in D : \alpha_2 w = \beta_2(w)=\eps(x), \beta_1 (w) =\eps(x) \} \]
and
\[ \Pi(D, V, x) = \{ v\in D : \alpha_1 (v) = \beta_1(v)=\eps(x), \alpha_2 (v) =\eps(x) \} \]
which have group structures induced from $+_2$, and $+_1$. Then
$\Pi(D, H) = \{\Pi(D, H, x):x\in X\}$ and $\Pi(D, V) = \{\Pi(D, V,
x):x\in X\}$  are totally  intransitive groupoids over $X$.

Clearly  maps
\[\delta_H : \Pi(D,H)\to H  \ \  \mbox{and} \ \ \delta_V: \Pi(D, V)\to V    \]
defined by $\delta_H (w) = \alpha_1 (w) $  and $\delta_V (v) =
\alpha_0(v),$ respectively, are morphisms of groupoids.

For a double groupoid  ${\cal D} =(D, H, V, X)$  it is shown in
\cite{Br-Sp1} that
\[   \gamma (\D) = (\Pi(D, H), H, \varepsilon )   \ \ \ \gamma' (\D) = (\Pi(D, V), V, \partial  )  \]
may be given the structure of crossed modules (see \cite{Br3} for
an exposition). So $\gamma$ is a functor from the category of
double groupoids to the category of crossed modules.

As we wrote, an   {\it edge symmetric double groupoid} is  a
double groupoid ${\cal D}$ but with the extra condition that the
horizontal and vertical edge groupoid structures coincide. These
double groupoids will, from now on, be our sole concern, and for
these it is convenient to denote the sets of points, edges and
squares by $X$, $G$, $D$. The identities in $G$ will be written
$\eps(x), 1_x$ or simply $1$. The source and target maps $G\to X$
will be written $\alpha$, $\beta$.

By a morphism $f:{\cal D}\to {\cal D'}$ of edge symmetric double
groupoid s is meant functions $f:D\to D'$, $f:G\to G'$, $f:X\to
X'$ which commute with all three groupoid structures.

A {\bf connection} for the double groupoid ${\cal D}$ is a
function $\Upsilon :G\to D$ such that if $a\in G$ then
$\Upsilon(a)$ has boundaries given by the following diagram $$
\spreaddiagramrows{2pc} \spreaddiagramcolumns{2pc}
\def \objectstyle {\ssize} \def \labelstyle{\sssize}
\def\labelstyle{\textstyle}\xymatrix{
\bullet \ar [d]_{a}  \ar @{}[dr] |{\Upsilon(a)} \ar [r]^{a} &
\bullet \ar [d]^{1}
\\ \bullet   \ar [r] _{1}                & \bullet
}  $$ and $\Upsilon$ satisfies the {\it transport law}:
  if $a,b\in G$ and $a+b$ is defined then
\begin{equation*}   \Upsilon (a+b)= (\Upsilon(a)+_1 \eps_2 b)+_2 \Upsilon(b)
\tag{*}\end{equation*} For further information on the transport
law and
 its uses, see \cite{Br-Mo}.

A morphism $f:{\cal D}\to {\cal D}'$ of edge symmetric double
groupoid  with special connections $\Upsilon, \Upsilon'$ is said
to preserve the connections if $f_2\Upsilon' = \Upsilon f_1$.

The category $DGrpd!$ has objects the pairs $({\cal D}, \Upsilon)$
of an edge symmetric double groupoid ${\cal D}$ with  connection
$\Upsilon$, and arrows the morphisms of edge symmetric double
groupoids preserving the connection. If $({\cal D}, \Upsilon)$ is
an object of $DGrpd!,$ then we have a crossed module $\gamma({\cal
D})$. Clearly $\gamma $ extends to a functor also written $\gamma
$ from  $DGrpd!$ to $CrsMod$, the category of crossed modules.

We now  show how  edge symmetric double groupoids arise from
crossed modules over groupoids.

Let ${\cal C}=(C, G, \delta)$ be a crossed module over groupoids
with base set  $X$. We define a edge symmetric double groupoid
${\cal D}({\cal C})$ as follows. First, $H=V=G$ with its groupoid
structure, base set $X$. The set  ${\cal D}({\cal C})$  of squares
is to consist of quintuples $$ w = \left(
{\spreaddiagramrows{-2pc}  \spreaddiagramcolumns{-1.5pc}
\objectmargin{0.05pc} w_1 : \diagram  & d & \\
                                     b & & c \\ & a &   \enddiagram }\right) \ \ \ \ \
{\spreaddiagramrows{2pc} \spreaddiagramcolumns{2pc}
\def \objectstyle {\ssize} \def \labelstyle{\sssize}
\def\labelstyle{\textstyle}\diagram
\bullet \dto^{b\ \ \ \ \ w}      \rto^{d}                &\bullet
\dto^{c} \\
  \bullet       \rto_{a}                &\bullet\todr^{w_1}
\enddiagram}
$$ such that $w_1\in C, a, b, c, d\in G$ and
\[            \delta (w_1) = -a-b+d+c.    \]
The vertical and horizontal structure on the set ${\cal D}({\cal
C})$ can be defined as follows.

The source and target maps on $w$ yield $d$ and $a$, respectively,
and vertical  composition is \begin{alignat*}{2}  \left(
{\spreaddiagramrows{-2pc} \spreaddiagramcolumns{-1.5pc}
\objectmargin{0.05pc} w_1 : \diagram  & d  & \\
                                     b&    &  c \\ & a &   \enddiagram }\right)&+_1
\left({ \spreaddiagramrows{-2pc}  \spreaddiagramcolumns{-1.5pc}
\objectmargin{0.05pc} w'_1 : \diagram  & a & \\
                                     b' & & c' \\ & a' &   \enddiagram }\right)&&=
\left({ \spreaddiagramrows{-2pc}  \spreaddiagramcolumns{-1.5pc}
\objectmargin{0.05pc} w'_1+{w_1}^{c'} : \diagram  & d & \\
                                   b+b' & & c+c' \\ & a' &   \enddiagram }\right),
                                   \\
\intertext{For the horizontal structure, the source and target
maps on $w$ yield $b$ and $c$, and the compositions are:} \left(
{\spreaddiagramrows{-2pc}  \spreaddiagramcolumns{-1.5pc}
\objectmargin{0.05pc} w_1 : \diagram  & d & \\
                                     b & &  c \\ & a &
                                     \enddiagram}
                                     \right)&+_2
\left({ \spreaddiagramrows{-2pc}  \spreaddiagramcolumns{-1.5pc}
\objectmargin{0.05pc} v_1 : \diagram  & e & \\
                                     c & & j \\ & i &   \enddiagram }\right)
&&=\left( {\spreaddiagramrows{-2pc}  \spreaddiagramcolumns{-1.5pc}
\objectmargin{0.05pc}{w_1}^i + v_1: \diagram  & d+e & \\
                                     b &    & j \\ & a+i &   \enddiagram} \right).
\end{alignat*}
Then ${\cal D}({\cal C})$ becomes a double groupoid with  these
structures. A connection on ${\cal D}({\cal C})$ is given by $$
\Upsilon(a)= \left({\spreaddiagramrows{-2pc}
\spreaddiagramcolumns{-1.5pc} \objectmargin{0.05pc} 1: \diagram  &
a & \\
                                     a &    & 1 \\ & 1 &
                                     \enddiagram}\right).
$$

The main result on double groupoids is:
\begin{thm}
The functor $\gamma :DGrpd!\to CrsMod$ is an equivalence of
categories \cite{Br-Sp1}.
\end{thm}

We emphasise that \cite{Br-Mac} contains a considerable
generalisation of this result.

We introduce a definition of linear coadmissible section for the
special double groupoid ${\cal D}({\cal C})$ as follows.
\begin{defn}
Let ${\cal C} =(C, G, \delta )$ be a crossed module and let ${\cal
D}({\cal C})$ be the corresponding  double groupoid. A {\bf linear
coadmissible section} $\sigma = (\sigma_0, \sigma_1):G\to {\cal
D}({\cal C})$ of ${\cal D}({\cal C})$ also written $$ \sigma(a) =
\left( \spreaddiagramrows{-2pc}  \spreaddiagramcolumns{-1.5pc}
\objectmargin{0.05pc}         \sigma_1(a) : \diagram     &
& \\
                                    \sigma_0\alpha(a)  &    &  \sigma_0\beta(a)\  \\ & a  &   \enddiagram \right)$$
is a pair of maps

\[  \sigma_0 :X\to G,  \ \ \ \ \sigma_1 :G\to  C \]
such that

\begin{enumerate}[(i)]
  \item if $x\in X$, $\beta \sigma_0(x) = x$, and if $a\in G$, then $\beta \sigma_1 (a) = \beta a$.
  \item if $a, b, a+b \in G$, then
\[    \sigma (a+b) = \sigma (a) +_2 \sigma (b) \]
  \item  $\alpha \sigma_0  :X\to X$ is a bijection,
$\alpha_1 \sigma  :G\to G$  is an automorphism.
\end{enumerate}

\end{defn}

Let ${\sf LinSec}({\cal D}({\cal C}))$ denotes the set of all
linear coadmissible sections. Then a group structure on ${\sf
LinSec}({\cal D}({\cal C}))$ is defined by the multiplication
\begin{align*}
 (\sigma * \tau)(z)& = \begin{cases} (\sigma_0\alpha\tau_0(z))+\tau_0(z) & \text {if } z\in
X,  \\   (\sigma\alpha_1(\tau(z)))+_1 \tau(z), & \text{if } z\in G
\end{cases}
\end{align*}
for $\sigma, \tau\in {\sf LinSec}({\cal D}({\cal C}))$.

It is easy to prove that the groups of linear coadmissible section
and free invertible derivation maps (coadmissible homotopies) are
isomorphic.


\section{Local coadmissible homotopies and \\ local linear sections }

Our aim in this section is to `localise' the concept of
coadmissible homotopy given in the previous section, analogously
to the way  Ehresmann \cite{Eh1} in the  $1$-dimensional case
localises a groupoid element to a local coadmissible section.

In order to cover both the topological and differentiable cases,
we use the term $C^r$ manifold for $r\geq -1$, where the case $r =
-1$ deals with the case of topological spaces and continuous maps,
with no local assumptions, while the case $r\geq 0$ deals as usual
with $C^r$ manifolds and $C^r$ maps. Of course, a $C^0$ map is
just a continuous map. We then abbreviate $C^r$ to  {\it smooth}.
The terms {\it Lie group } or {\it Lie groupoid} will then involve
smoothness in this extended sense.  By a {\it partial
diffeomorphism } $f:M\to N$ on $C^r$ manifolds $M, N$ we mean an
injective partial function with open domain and range and such
that $f$ and $f^{-1}$ are smooth.

One of the key differences between the cases $r = -1$ or $ 0$ and
$ r\geq 1$ is that for $ r\geq 1$, the pullback of $C^r$ maps need
not  be a smooth submanifold of the  product, and so
differentiability of maps on the pullback cannot always be
defined. We therefore adopt the following definition of Lie
groupoid. Mackenzie \cite{Mac1} discusses the utility of various
definitions of differential groupoid.

A {\bf Lie groupoid} is  a topological groupoid $G$ such that
\begin{enumerate}[(i)]
  \item the space of arrows is a smooth manifold, and the space of objects
is a smooth submanifold of $G$,

  \item  the source and target maps $\alpha, \beta$ are smooth maps
and are submersions.
  \item  the domain $G\sqcap_{\beta } G$ of the difference map is a smooth submanifold
of $G\times G$, and
\item the difference map $ \partial: (a,b) \mapsto a-b$ is a smooth map.
\end{enumerate}

Recall that coadmissible homotopies were defined in the previous
section. Here we define the local version.

\begin{defn}
 Let ${\cal C} = (C, G, \delta)$ be a  crossed module such that
$(G, X)$ is a Lie groupoid. A {\bf local coadmissible homotopy }
$s = (s_0, s_1)$ on $U_0, U_1$ consists of two partial maps
\[  s_0: X \to G  \  \ \  \ s_1:G\to C        \]
with open domains $U_0\subseteq X$, $U_1\subseteq G$, say, such
that $\alpha (U_1),$  $\beta (U_1)\subseteq U_0$ and
\begin{enumerate}[(i)]
  \item If $x\in U_0$, then  $\beta s_0(x) =x$.
  \item If $a, b, a+b\in U_1$, then
\[   s_1(a+b)  = s_1(a)^b+s_1(b),      \]
we say $s_1$ is a local derivation.
  \item If $a\in U_1$ then $\beta s_1(a) = \beta (a)$.
  \item If $\Delta(s)=(f_0, f_1,f_2)$ is defined by \begin{align*}
  f_0 (x) &= \alpha s_0 (x), \ \ x\in U_0, \\
f_1 (a) &= s_0 \alpha (a) + a + \delta s_1 (a) -s_0\beta (a),a\in
U_1\\ f_2(c)&= (c+\delta s_1(c))^{-s_0\beta(c)}, c \in C
\end{align*}
then $f_0, f_1$ are partial diffeomorphisms  and ${f_1}^{-1}$,
$f_1$ are linear.
\end{enumerate}
\end{defn}
Note that $f_1$ injective implies $f_2$ injective by Theorem
\ref{14}. We cannot put smooth conditions on $f_2$ as we do not
have a topology on $C$.

A local coadmissible homotopy  $s$ defined as above can be
illustrated by the following diagram. $$ \spreaddiagramrows{2pc}
\spreaddiagramcolumns{2pc}
\def\labelstyle{\textstyle}\diagram
         C \dto_{\delta} \rto ^{f_2} & C \dto^{\delta} \\
        U_1  \urto|{s_1}\dto \rto ^{f_1}     & G \dto              \\
        U_0  \urto|{s_0}          \rto _{f_0} & X
\enddiagram
$$

Suppose given  open subsets $V_0\subseteq X$ and  $V_1\subseteq G$
such that $\alpha(V_1), \beta(V_1)\subseteq V_0$. Let $t:(V_0,
V_1)\to (C, G)$ be a local coadmissible homotopy on $V_0$, $V_1$
with $\Delta(t) =g$. Let $\Delta(s)=f$ be as above. Now we can
define a multiplication of $s$ and $t$ in the following way $$ (s
* t)(z)= \begin{cases}t_1(z) + s_1g_1(z)^{t_0(\beta z)},  &   z \in U_1
\\
 s_0 g_0(z) + t_0(z), & z\in U_0.
\end{cases}$$
We write $\Do(s_\eps)$ for the domain of a function $s_\eps$.
\begin{lem}
The product  $s * t$ is a local  coadmissible homotopy.
\end{lem}
\begin{pf}
We will prove that the domain of $s *t$ is open. In fact, if $a\in
U_1$, $g_1(a)\in V_1$, $\beta (a) \in U_0$, then $a\in U_1\cap
g_1^{-1}(V_1)\cap \beta^{-1}(U_0)$ is an open set in $G$ and also
if $x\in U_0$ and $g_0(x)\in V_0$ then $x\in U_0\cap
g_0^{-1}(V_0)$ is open in $X$, so the domain of $(s *t)$ is open.
One can show that $$\beta (s *t)_0(x) = \beta (x),   \ \mbox{for }
x\in \Do(s*t)_0$$ $$\beta (s *t)_1(a) = \beta (a), \  \mbox{for }
a\in \Do(s*t)_1$$ and $(s * t)_1$ is a derivation map as in
Proposition 2.4 in \cite{Br-Ic2}. i.e.,
\[  (s * t)_1(a+b) = (s * t)_1(a)^b + (s * t)_1(b)     \]
for $a, b, a+b\in \Do(s*t)_1$. We define maps $h_0, h_1$ as
follows:
\[ h_0(x) = f_0g_0(x)= \alpha (s*t)_0(x) \ \ \  \mbox{for} \ \  x\in \Do(s*t)_0 \]
\[ h_1(a) = f_1g_1(a)= (s*t)_0(\alpha a) +a+\delta(s*t)_1(a)-(s*t)_0(\beta a) , \ \  \mbox{for } a\in \Do(s*t)_1.\]
Since $h_0, h_1$ are compositions of partial diffeomorphisms, they
are partial diffeomorphisms.
\end{pf}

\begin{prop}
Let ${\sf LFDer}^*({\cal C} )$ denotes the set of all local
coadmissible homotopies of a crossed module ${\cal C} = (C, G,
\delta)$ such that $(G, X)$ is a Lie groupoid. For each $s, t\in
{\sf LFDer}^*({\cal C}) $, $s * t \in {\sf LFDer}^*({\cal C}) $
and for each $s\in {\sf LFDer}^*({\cal C}) $, let $\Delta(s)=f$
and let
\begin{equation*}
s^{-1}(z)=  \begin{cases} -s_1({f_1}^{-1}(z))^{{s_0}^{-1}(\beta
z)} & \text{if } \ \ z\in U_1, \\
 -s_0({f_0}^{-1}(z))  &  \text{if} \ \  z\in U_0.
\end{cases}
\end{equation*}
Then $s^{-1}\in {\sf LFDer}^*({\cal C}),$ and with this product
and  inverse element the set ${\sf LFDer}^*({\cal C}) $ of local
coadmissible homotopies becomes an inverse semigroup.
\end{prop}
\begin{pf} The proof is fairly straightforward.
\end{pf}

Recall that linear coadmissible sections  were defined in the
previous section. Here we define the local version.
\begin{defn}{\rm
Let ${\cal C} = (C, G,\delta)$ be a crossed module of groupoids
with base space $X$ and let ${\cal D}({\cal C})$ be the
corresponding edge symmetric double groupoid such that $(G, X)$ is
a Lie groupoid.

A {\bf local linear section} $\sigma = (\sigma_0, \sigma_1):G\to
{\cal D}({\cal C}),$ written $$ \sigma(a) =  \left(
\spreaddiagramrows{-2pc}  \spreaddiagramcolumns{-1.5pc}
\objectmargin{0.05pc} \sigma_1(a) : \diagram     & &
\\ \sigma_0\alpha(a) &    &  \sigma_0\beta(a)\  \\ & a  &
\enddiagram \right)$$ consists of two partial maps
\[  \sigma_0: X \to G, \qquad  \sigma_1:G\to  C        \]
with open domains $U_0\subseteq X$, $U_1\subseteq G$, say, such
that $\alpha (U_1),$  $\beta (U_1)\subseteq U_0$ and
\begin{enumerate}[(i)]
  \item if $x\in U_0$, then  $\beta \sigma_0(x) =x$, and if $a\in U_1$,
then  $\beta \sigma_1 (a) = \beta a$;
  \item if $a, b, a+b\in U_1$, then
\[   \sigma(a+b)  = \sigma(a)+_2\sigma(b)      \]
(this is the local linear condition);
  \item if $f_0, f_1$ are defined by \begin{align*} f_0 (x)& = \alpha \sigma_0 (x), \ \ x\in
  U_0, \\ f_1 (a) & = \alpha_1 \sigma (a) , \ a\in U_1,
  \end{align*} then $f_0, f_1$ are partial diffeomorphisms and $f_1, {f_1}^{-1}$
are linear.
\end{enumerate}
 }
\end{defn}

Given open subsets $V_0\subseteq X$ and  $V_1\subseteq G$ such
that $\alpha(V_1), \beta(V_1)\subseteq V_0$, let $\tau $ be a
local linear section with domain $V_0$ and $V_1$. Let $\sigma$ be
as above. Now we can define a multiplication of $\sigma$  and
$\tau $ in the following way $$ (\sigma *\tau )(z)= \begin{cases}
\sigma(\alpha_1 \tau)(z)+_1 \tau (z),  & z\in U_1 \\ \sigma_0
(\alpha \tau_0)(z) + \tau_0 (z), & z\in U_0
\end{cases}$$

\begin{lem}

The product function $\sigma *\tau$  is a local linear
coadmissible section.
\end{lem}
\begin{pf}
The key point is to prove that the  domain of $\sigma *\tau$  is
open. In fact, if $a\in U_1$, $\alpha_1 \sigma(a)\in V_1$, $\beta
(a) \in U_0$, then $a\in U_1\cap (\alpha_1 \sigma)^{-1}(V_1)\cap
\beta^{-1}(U_0)$ is an open set in $G$ and also if $x\in U_0$ and
$\alpha \sigma_0(x)\in V_0$ then $x\in U_0\cap (\alpha
\sigma_0)^{-1}(V_0)$ is open in $X$, so domain of $(\sigma *\tau)$
is open.

The remaining part is easily done.\end{pf}
\begin{prop}
Suppose ${\cal C}=(C,G, \delta)$ is a crossed module such that $G$
is a Lie groupoid on $X$. Let $\Gamma({\cal D}({\cal C}))$ denote
the set of all local linear coadmissible sections of ${\cal
D}({\cal C})$. For each $\sigma, \tau\in \Gamma ({\cal D}({\cal
C})) $, $\sigma * \tau \in \Gamma ({\cal D}({\cal C})) $ and for
each $\sigma\in \Gamma ({\cal D}({\cal C})) $, let

\begin{equation}
\sigma^{-1}(z)=\left \{ \begin{array}{ll}
 -_1\sigma({\alpha \sigma})^{-1}(z) &  \text{if} \  \ z\in U_1, \\
-\sigma_0({\alpha \sigma_0})^{-1}(z))  &\text{if}  \ \ z\in U_0.
\end{array}
\right.
\end{equation}
Then  with this product and  inverse element, the set
$\Gamma({\cal D}({\cal C})) $ of local linear coadmissible
sections becomes an inverse semigroup.
\end{prop}
\begin{pf}
The inverse  of $\sigma$ is a linear coadmissible section, because
it is a composition of the linear maps $\sigma$ and $(\alpha
\sigma)^{-1}$.
\end{pf}
\begin{prop}\label{27}
Let $\Delta(s_0, s_1)=f$ be a local coadmissible homotopy for a
crossed module ${\cal C} =(C, G, \delta)$ and let ${\cal D}({\cal
C})$ be the corresponding double groupoid. A partial map  $s$ is
defined by
\[       s= (s_0, s_1): G\to {\cal D}({\cal C})   \]
$$  s(a)= (s_0, s_1)(a) = \left( \spreaddiagramrows{-2pc}
\spreaddiagramcolumns{-1.5pc}
   \objectmargin{0.05pc} s_1(a) : \diagram  & f_1(a) & \\
                                   s_0(x) &         & s_0(y) \\ & a &   \enddiagram \right).$$
Then $s$  is a local linear coadmissible section.
\end{prop}
\begin{pf}
It is easy to see from the definitions  of local coadmissible
sections, local coadmissible homotopies \cite{Br-Ic2}.
\end{pf}
\begin{cor}
The inverse semigroups of local coadmissible homotopies and local
linear sections are isomorphic.
\end{cor}

Throughout the next sections, we will deal with linear
coadmissible sections rather than  coadmissible homotopies.

\section{V-locally Lie double groupoids}

Let ${\cal C} = (C, G, \delta)$ be a  crossed module such that
$(G, X)$ is a Lie groupoid. Let ${\cal D}({\cal C})$ be the
corresponding double groupoid. Let $\Gamma ({\cal D}({\cal C}))$
be the set of all local linear coadmissible sections and let $W$
be a subset of $C$ such that $W$  has the structure of a manifold
and $\beta :W\to X$ is a smooth surmersion.

Let  $W^G$ be the set \begin{equation*} \{
  \left( {\spreaddiagramrows{-2pc}  \spreaddiagramcolumns{-1.5pc}
\objectmargin{0.05pc} w_1 : \diagram  & d & \\
                                     b & & c \\ & a &   \enddiagram} \right)
: w_1\in W, a, b, c\in G,\beta (b) = \alpha (a), \beta (a) = \beta
(c) = \beta (w_1), d=b+a+\delta(w_1)-c  \}  .\tag{*}
\end{equation*}  Here the set $W^G$ can be considered as
a repeated pullback, i.e., if
\[ G^3 =\{ (b, a, c) : \alpha(a)= \beta(b), \beta(a)=\beta(c) \}  \]
is a pullback, then $$ \diagram
        W^G\dto \rto  & G^3 \dto^{\beta \pi_1} \\
        W \rto_{\beta} & X
\enddiagram
$$ is a pullback, so $W^G$ has a manifold structure on it, because
$\beta$ and $\beta \pi_1$ are smooth and surmersions.

We can show each element $w\in W^G$ by the following diagram: $$
\spreaddiagramrows{2pc}
\spreaddiagramcolumns{6pc}\objectmargin{0pc}
\def\objectstyle{\ssize} \def\labelstyle{\ssize}
\def\labelstyle{\textstyle}\xymatrix{
\bullet \ar [d]^{b} \ar @{} [dr] |w \ar @{.>}
[r]^{b+a+\delta(w_1)-c} & \bullet \dto^{c}   \\ \bullet^{} \ar [r]
_{a} & \bullet \ar @(d,r)[]_{w_1} } $$ Clearly $W^G\subseteq
W\times G\times G\times G$ and $W^G\subseteq {\cal D}({\cal C})$.

A local linear  coadmissible section  $ s=(s_0, s_1)$ as given in
Proposition \ref{27}  is said to be smooth if $Im (s_1)\subseteq
W$ and $s_0, s_1$ are smooth.  We emphasise that in our current
context, such an $s$ should be regarded as a `local procedure'.
Let $\Gamma^r(W^G)$ be the set of local linear smooth coadmissible
sections. We say that the triple $(\alpha_1, \beta_1, W^G) $ has
{\bf enough smooth local linear coadmissible sections} if for each
$ w = \left( \spreaddiagramrows{-2pc}
\spreaddiagramcolumns{-1.5pc} \objectmargin{0.05pc} w_1 : \diagram
& d & \\
                                     b & & c \\ & a &
\enddiagram \right)\in W^G$, there is a local linear  smooth
coadmissible section  $\Delta(s)=f$ with domains  $(U_0,U_1)$ such
that

\begin{enumerate}[(i)]
  \item   $s\beta_1(w) = w$, $\alpha_1(w) = f_1(a)=d;$  $s_1 \beta_1(w) = w_1 = s_1 (a)$,
$s_0 \beta (a) = c, s_0 \alpha (a) = b$.
  \item  the values of $s $ lie in $W^G$
  \item $s$ is  smooth as a pair of function
$U_0 = \Do(s_0)\to G$ and $U_1 = \Do(s_1)\to W^G$.
\end{enumerate}
We call such an $s$ a {\bf local linear smooth coadmissible
section through $w$.}

\begin{defn}\label{llg}
Let  ${\cal C} = (C, G, \delta)$ be a crossed module  over a
groupoid with base space $X$ and let  ${\cal D}({\cal C})$ be the
corresponding double groupoid. A {\bf $V$-locally Lie double
groupoid} structure $({\cal D}({\cal C}), W^G)$ on  ${\cal
D}({\cal C})$ consists of a smooth structure on $G, X$ making $(G,
X)$ a Lie groupoid and a smooth manifold $W$ contained as a set in
$C$ such that if $W^G$ is as in (*), then
\begin{enumerate}[S1)]
  \item  $W^G = -_1{W^G}$;
  \item  $G\subseteq W^G\subseteq {\cal D}({\cal C})$;
  \item  the set $(W^G\sqcap_{\beta_1} W^G)\cap {\sf d}^{-1}(W^G) = {W_{\sf d}}^G$
is open in $(W^G\sqcap_{\beta} W^G)$    and the restriction to
${{W_{\sf d}}}^G$ of the difference map
\[   {\sf d} : {\cal D}({\cal C})\sqcap_{\beta_1} {\cal D}({\cal C})\to {\cal D}({\cal C})     \]
\[     (w, v)\mapsto w -_1 v, \]
is smooth;
\item the restriction to $W^G$ of the source and target maps
$\alpha_1$
 and
$\beta_1 $ are smooth and the triple $(\alpha_1, \beta_1, W^G)$
has enough local linear smooth coadmissible sections,
\item $W^G$ generates ${\cal D}({\cal C})$ as a groupoid with respect to $+_1$.
\end{enumerate}

\end{defn}

Also one can define locally Lie crossed module structure on a
crossed module by considering the above Definition \ref{llg}.

Let ${\cal C} = (C, G, \delta)$ be a crossed module over a
groupoid with base space $X$. A {\bf locally Lie crossed module
structure}  $(C, W, \delta)$ on ${\cal C}$ consists of a smoooth
structure on $G$, $X$ making $(G, X)$  a Lie groupoid and a subset
$W$ of $C$ with a smooth structure on $W$ such that $W$ is
$G$-equivariant and
\begin{enumerate}[C1)]

\item  $(C, W)$ is a locally Lie groupoid,

\item  $I(G)\subseteq W\subseteq C$,

\item  the restriction to $W$ of the map $\delta:C\to G$ is smooth,

\item  the set $W_A= A^{-1}(W)\cap(W\sqcap_{\beta} G)$ is open in $W\sqcap_{\beta}G$ and the
restriction to $A_W: W_A\to W$ of the action $A:C\sqcap_{\beta}
G\to C$ is smooth.

\item  Let

\[ \underline W = \{ (w;b,a,c): w\in W, a,b,c\in G, \beta(a) =\beta(w), \alpha(a) = \beta(b), \beta(c) = \beta(a) \}.   \]

We say that $\underline W$ has enough local smooth coadmissible
homotopies if for all $(w; b, a, c)\in {\underline W}$ there
exists a local smooth coadmissible homotopy $(s_0, s_1)$ such that
$s_1(a) = w, s_0\beta (a) = c, s_0\alpha (a) = b.$
\end{enumerate}
Let us compare the above two definitions.

First of all, in the definition of locally Lie crossed module,
conditions $C3$   and $C4$ gives rise to the difference map
\[   {\sf d} : {\cal D}({\cal C})\sqcap_{\beta_1} {\cal D}({\cal C})\to {\cal D}({\cal C})     \]
\[     (w, v)\mapsto w -_1 v, \]
which is smooth. In fact,

\begin{align*} {{\sf d} \left( { \left( { \spreaddiagramrows{-2pc}  \spreaddiagramcolumns{-1.5pc}
                                                  \objectmargin{0.05pc} w_1 : \diagram     & d  & \\
                                                  b &    &  c  \\ & a &   \enddiagram }\right)} ,
                                                  {  \left( \spreaddiagramrows{-2pc}  \spreaddiagramcolumns{-1.5pc}
                                                  \objectmargin{0.05pc} v_1 : \diagram     & d'  & \\
                                                    b' &    &  c'  \\ & a &   \enddiagram \right)}\right) }
                                                  & = {  \left( \spreaddiagramrows{-2pc}  \spreaddiagramcolumns{-1.5pc}
                                                  \objectmargin{0.05pc} w_1 : \diagram     & d  & \\
                                                  b &    &  c  \\ & a &   \enddiagram \right)}-_1
                                                  {   \left( \spreaddiagramrows{-2pc}  \spreaddiagramcolumns{-1.5pc}
                                                  \objectmargin{0.05pc} v_1 : \diagram     & d' & \\
                                                    b' &    &  c'  \\ & a &   \enddiagram \right) } \\
                                                  & = {  \left( \spreaddiagramrows{-2pc}  \spreaddiagramcolumns{-1.5pc}
                                                  \objectmargin{0.05pc} {w_1} : \diagram     & d  & \\
                                                  b &    &  c  \\ & a &   \enddiagram \right)}+_1
                                                  {   \left( \spreaddiagramrows{-2pc}  \spreaddiagramcolumns{-1.5pc}
                                                  \objectmargin{0.05pc} {-v_1}^{-c'} : \diagram     & a  & \\
                                                    -b' &    &  -c'  \\ & d' &   \enddiagram
                                                    \right)}
\\& =   \left({ \spreaddiagramrows{-2pc}  \spreaddiagramcolumns{-1.5pc}
                                       \objectmargin{0.05pc} ({-v_1}')^{-c'}+{w_1}^{-c'} : \diagram     & d  & \\
                                                  b-b' &    &  c-c'  \\ & d' &   \enddiagram } \right)
                                                 \end{align*}

Since $\delta_W$, $+$, $A_W$ are smooth, ${\sf d}$ is smooth. This
is equivalent to the two smooth conditions $C3$, $C4$ for locally
Lie crossed module, because the  formulae for ${\sf d}$ involves
$+$ and the action $A_W$.

The condition $C1$ that $(C, W)$ is a locally Lie groupoid, which
includes $W$ generates $C$. The other equivalent condition can be
stated as follows: We first prove that if $W$ generates $C$ and is
$G$-equivariant, then $W_G$ generates ${\cal D}({\cal C})$ with
respect to $+_1$.

Let $ w = \left( \spreaddiagramrows{-2pc}
\spreaddiagramcolumns{-1.5pc}
   \objectmargin{0.05pc} \gamma : \diagram  & d & \\
                                   b &         & c \\ & a &   \enddiagram \right)\in {\cal D}({\cal C})$.
We prove by induction that if $\gamma$ can be expressed as a word
of length $n$ in conjugates of elements of $W$ then $w$ can be
expressed as
\[       w = w^1 +_1\cdots+_1w^n   \]
where
 $ w^i = \left( \spreaddiagramrows{-2pc}  \spreaddiagramcolumns{-1.5pc}
   \objectmargin{0.05pc} \gamma_i : \diagram  & d_i & \\
                                   b_i &         & c_i \\ & a_i &   \enddiagram \right)\in W^G,$
for $i=1,\ldots ,n$, and $\gamma_i\in W$.

This is certainly true for $n = 1$, since $w\in W^G$ if and only
if $\gamma \in W$.

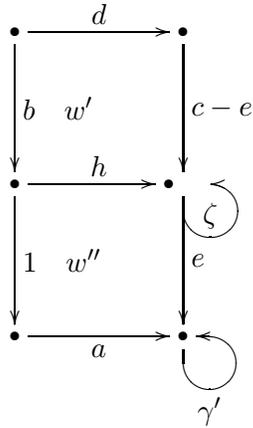
$$ \spreaddiagramrows{2pc} \spreaddiagramcolumns{2pc}
\def \objectstyle {\ssize} \def \labelstyle{\sssize}
\def\labelstyle{\textstyle}\diagram
\bullet \dto^{b\ \  \ w'} \rto^d            & \bullet \dto^{c-e}
\\ \bullet \dto^{1\ \ \  w''}  \rto^{h}          & \bullet
\todr^{\zeta} \dto^{e}     \ \ \       \\ \bullet     \rto_{a}
& \bullet \todr_{\gamma'}
\enddiagram
$$ Suppose $\gamma = \gamma'+\zeta^e$ where $\gamma'$ can be
expressed as a word of length  $n$ in conjugates of elements of
$W$ and $\zeta \in W$.

Let $h = a+\delta \gamma'-e,$ and let $  w'' = \left(
\spreaddiagramrows{-2pc}  \spreaddiagramcolumns{-1.5pc}
   \objectmargin{0.05pc} \gamma' : \diagram  & h & \\
                                   1 &         & e \\ & a &   \enddiagram \right).$
Then $w''\in {\cal D}({\cal C})$ and so $w''$ can be expressed as
a word of length $n$ in elements of $W^G$, by the inductive
assumption.

Let $ w' = \left( \spreaddiagramrows{-2pc}
\spreaddiagramcolumns{-1.5pc}
   \objectmargin{0.05pc} \zeta : \diagram  & d & \\
                                   b &         & c-e \\ & h &   \enddiagram \right).$
Then $w'\in {\cal D}({\cal C}),$ since
\begin{align*}
\delta \zeta & =  \delta (-\gamma'+\gamma)^{-e}\\
               & =  -\delta (\gamma')^{-e}+\delta (\gamma)^{-e}\\
               & =  e-\delta \gamma'-e+e+\delta \gamma -e\\
               & =  e-(-a+h+e) -a -b+d+c -c\\
               & =  e-e-h+a-a-b+d+c-e\\
               & =  -h -b+d+c-e
\end{align*}
and $\zeta\in W$.  Clearly $w = w'+_1 w'',$ and so $w$ can be
expressed length $n+1$ in conjugates of elements of $W^G$.

Conversely, suppose  $W^G$ generates ${\cal D}({\cal C})$ with
respect to $+_1$.

Let $\gamma\in C$ and let $ w = \left( \spreaddiagramrows{-2pc}
\spreaddiagramcolumns{-1.5pc}
   \objectmargin{0.05pc} \gamma : \diagram  & 1 & \\
                                   1 &         & 1 \\ & 1 &   \enddiagram \right).$
Then $w\in {\cal D}({\cal C)}$. Since $W^G$ generates ${\cal
D}({\cal C})$, we can write
\[       w = w^1 +_1\cdots+_1w^n   \]
where $  w^i = \left( \spreaddiagramrows{-2pc}
\spreaddiagramcolumns{-1.5pc}
   \objectmargin{0.05pc} \gamma_i : \diagram  &  & \\
               b_i &         & c_i \\ & a_i &   \enddiagram \right), $
$\gamma_i \in W$, for $i =1,\ldots, n$ and $w^i\in W_G.$ Then $$ w
= \left( \spreaddiagramrows{-2pc}  \spreaddiagramcolumns{-1.5pc}
   \objectmargin{0.05pc} {\gamma_1}^{c_1}+{\gamma_2}^{c_2}+\cdots+{\gamma_n}^{c_n+\cdots+c_1} : \diagram  & d & \\
                             b_1+\cdots+b_n       &         & c_1+\cdots+c_n \\ & a &   \enddiagram \right) =
   \left( \spreaddiagramrows{-2pc}  \spreaddiagramcolumns{-1.5pc}
   \objectmargin{0.05pc} \gamma : \diagram  & 1 & \\
                                   1 &         & 1 \\ & 1 &   \enddiagram \right)$$
We get
\[  \gamma = {\gamma_1}^{c_1}+{\gamma_2}^{c_1+c_2}+\cdots+{\gamma_n}^{c_n+\cdots+c_1}, \ \ \gamma_i\in W  \]
So $W$ generates $C$.

In the definition of $V$-locally Lie double groupoid, condition S4
transfers as follows: Let $(\alpha_1, \beta_1, W^G)$ have enough
local linear smooth coadmissible sections. Then for each $ w =
\left( \spreaddiagramrows{-2pc}  \spreaddiagramcolumns{-1.5pc}
   \objectmargin{0.05pc} \gamma : \diagram  &  & \\
                                   b &         & c \\ & a &   \enddiagram \right)\in W^G$
there exists  a local linear smooth coadmissible section $s$ such
that $s\beta_1 (w) = w$, i.e., there exists $(s_0, s_1)$ that is a
local coadmissible homotopy for the crossed module ${\cal C} = (C,
G, \delta)$. So for $(w; b, a, c)$ defined as above, there exists
a local smooth coadmissible homotopy $s=(s_0, s_1)$ such that
$s_1(a) =w, s_0\alpha a = b, s_0\beta a = c.$

We state  some deductions from the axioms S1--S5.
\begin{enumerate}[1)]
\item The inverse map  ${\sf i}: W^G\to W^G, w \mapsto -_1w$ is smooth,
\item  The set ${\sf d}^{-1}(W^G)\cap({W^G}_{\alpha_1}\sqcap_{\beta_1}W^G)$ is open
${W^G}_{\alpha_1}\sqcap_{\beta_1}W^G$ and the horizontal product
map is smooth on this set,
\item Let $k, s, t\in \Gamma(W^G)$, and suppose
$a\in \Do((k*s*t)_1)$ and $b=\beta_1 t(a)$ satisfy $(k*s)(b)\in
W^G$ and $(k*s*t)(a)\in W^G$. Then there are restriction $k', s',
t'$ of $k, s, t$ such that $a\in \Do((k'*s'*t')_1)$ and
$k'*s'*t'\in \Gamma (W^G)$ (compare 1.6 in \cite{Aof-Br}).
\end{enumerate}
There is a lemma which we shall use later.
\begin{lem} \label{st}
Suppose $s, t \in \Gamma^r (W^G)$, $a\in G$ and $s (a) = t (a)$.
Then there is a pair of neighbourhoods $(U_0, U_1)$, where $U_0$
is a neighbourhood both of $\alpha (a)$ and $\beta (a)$ and $U_1$
is a neighbourhood of $a$ such that the restriction of $ s *
t^{-1}$ to $(U_0, U_1)$ lies in $\Gamma^r (W^G)$.
\end{lem}
\begin{pf}
Since $s$ and $t$ are smooth and $s(a) = t(a)$, then $(s(a),
t(a))\in W^G\sqcap_{\beta_1}W^G$. This gives rise to maps
\[  (s_0, t_0) : \Do(s_0 )\cap \Do(t_0 )\to G\sqcap_{\beta} G  \ \ \  \mbox{and}
 \ \ \ (s, t) : \Do(s_1 )\cap \Do(t_1 )\to W^G\sqcap_{\beta_1} W^G   \]
which are smooth. But by condition S of Definition \ref{llg},
$(W^G\sqcap_{\beta_1} W^G)\cap {\sf d}^{-1}(W^G)$ is open in
$W^G\sqcap_{\beta_1} W^G$ and $(G,X)$ is a (globally) Lie
groupoid.  Hence there exist  open neighbourhoods $U_1$ of $a$ in
$G$ , $U_0$ of $\alpha (a)$, $\beta (a) $ in $X$ such that $$(s,
t)(U_1)\subseteq (W^G\sqcap_{\beta_1} W^G)\cap {\sf d}^{-1}(W^G) \
\ \  \mbox{and}
 \ \ \ (s_0, t_0)(U_0)\subseteq (G\sqcap_{\beta} G)\cap {\partial}^{-1}(G).$$
Hence ${\sf d}(s, t)(U_1)$ is contained in $W^G$ and ${\partial}
(s_0, t_0) (U_0)$ is contained in $G$. This gives $({s } *t^{-1})
(U_1)\subseteq W^G$ and $({s } *t^{-1})_0 (U_0)\subseteq G$. So
${s } * t^{-1} \in \Gamma^r (W^G)$.
\end{pf}

\noindent {\bf \large Germs}

Let $s$, $t$  be two local linear   coadmissible sections with
domains  $( U_0, U_1)$ and $(U'_0,U'_1)$ respectively, and let
$a\in U_1\cap U'_1$. We will define an equivalence relation as
follows: set $s \sim_a t$ if and only if $U_1\cap U'_1$ contains
an  open neighbourhood $V_1$ of $a$ such that
\[        s \mid_{ V_1} =  t \mid_{ V_1}        \]
and $\alpha (V_1), \beta (V_1)\subseteq V_0.$

Let $J_a({\cal D}({\cal C}))$ be the set of all equivalence
classes of $\sim_{a}$ and let
\[        J({\cal D}({\cal C})) = \bigcup \{ J_a({\cal D}({\cal C})) : a\in G \}.       \]
Each element of $J_a({\cal D}({\cal C}))$ is called a {\bf germ}
at $a$ and is denoted by $[s]_{a}$ for $s\in \Gamma ({\cal
D}({\cal C})) $, and $J({\cal D}({\cal C}))$ is called the sheaf
of germs of local linear  coadmissible sections of the double
groupoid ${\cal D}({\cal C})$.
\begin{prop}
Let $J({\cal D}({\cal C}))$ denote the set of all germs of local
linear  coadmissible sections of the double groupoid ${\cal
D}({\cal C})$. Then $J({\cal D}({\cal C}))$ has a  natural
groupoid structure  over $G$.
\end{prop}
\begin{pf}
Let $s, t\in \Gamma ({\cal D}({\cal C})) $ and $ \Delta(s)= f$,
$\Delta(t)=g$. The source and target maps are \begin{align*}
\alpha ([s ]_{a}) &= f_1 (a)   \\    \beta ([s ]_{a} )  &= a    \\
\intertext{and the object map is $a\mapsto [1]_{a}$, the
multiplication $*$ is}     [s ]_{g_1a}*   [t ]_{a} &= [(s *t)]_{a}
\\ \intertext{and the inversion map is}
    [s ]^{-1}_{a} &= [s^{-1}]_{f_1a}.
\end{align*}
\end{pf}

\noindent {\bf Remark}: One can give  a sheaf topology on $J({\cal
D}({\cal C}))$ defined by taking as basis the sets $\{ [s ]_{a} :
a\in U_1 \}$ for $s\in \Gamma ({\cal D}({\cal C})) $, $U_1$ open
in $G$. With this topology  $J({\cal D}({\cal C}))$ is a
topological groupoid. We do not use the sheaf topology since this
will not give $W^G$ embedded as an open set.

Suppose now that $({\cal D}({\cal C}), W^G)$ is a $V$-locally Lie
double groupoid. Let $\Gamma^r(W^G)$ be the subinverse semigroup
of $\Gamma ({\cal D}({\cal C}))$ consisting of local  linear
coadmissible sections with values in $W^G$ and which are smooth.
Let $\Gamma^r({\cal D}({\cal C}), W^G)$ be the subinverse
semigroup of $\Gamma ({\cal D}({\cal C}))$ generated by
$\Gamma^r(W^G)$.

\noindent {\bf Remark:}  In view of the discussion in the
Introduction, we regard an element of this inverse semigroup as an
`iteration of local procedures'.

If $s \in \Gamma^r({\cal D}({\cal C}), W^G)$, then there are $s^i
\in \Gamma^r(W^G)$, $i =1,\cdots,n$ such that
\[     s =s^n*\cdots*s^1.   \]
So let $J^r({\cal D}({\cal C}))$ be the subsheaf of $J({\cal
D}({\cal C}))$ of germs of elements of $\Gamma^r({\cal D}({\cal
C}), W^G)$.  Then $J^r({\cal D}({\cal C}))$ is generated as a
subgroupoid of $J({\cal D}({\cal C}))$ by the sheaf $J^r(W^G)$ of
germs of element of $\Gamma^r(W^G)$. Thus an element of $J^r({\cal
D}({\cal C}))$ is of the form
\[    [s ]_a = [s^n]_{a_n} *  \cdots * [s^1]_{a_1}   \]
where   $s =s^n*\cdots*s^1$ with $[s^i]_{a_i}\in J^r(W^G)$,
$a_{i+1} = f_i (a_i), i = 1, \ldots,  n$ and $a_1 = a \in
\Do(s^1).$

Let $\psi : J({\cal D}({\cal C}))\to {\cal D}({\cal C})$ be the
final map defined by $$\psi ([s]_a) = s(a) = \left(
\spreaddiagramrows{-2pc}  \spreaddiagramcolumns{-1.5pc}
   \objectmargin{0.05pc} s_1(a) : \diagram  & f_1(a) & \\
                                   s_0(x) &         & s_0(y) \\ & a &   \enddiagram \right),$$
where $s$ is a local linear coadmissible section. Then  $\psi$ is
a groupoid morphism. In fact, let $ \Delta(s)=f$, $ \Delta(t)=g$,
then
\begin{align*}
\psi ([s]_{g_1(a)} * [t]_a)& = \psi ([s*t]_a)\\
                         & =  (s * t) (a) \\
                         & =  s\alpha_1 t(a) +_1 t(a)\\
                         & = s(g_1(a)) +_1 t(a) \\
                         & =  \psi [s]_{g_1(a)} +_1 \psi [t]_a.
\end{align*}
Then
\[       \psi (J^r({\cal D}({\cal C}))) = {\cal D}({\cal C}),   \]
from the axiom S4 of a $V$-locally Lie double groupoid on  ${\cal
D}({\cal C})$ in Definition \ref{llg}.

Let $J_0 = J^r(W^G) \cap Ker \ \psi$, where as usual
\[    Ker \ \psi = \{ [s ]_a : \psi [s ]_a = 1_a  \}   \]
We will prove that $J_0$ is a normal subgroupoid of $J^r({\cal
D}({\cal C}))$. This allows us to  define a quotient groupoid
$J^r({\cal D}({\cal C}))/J_0$ in the next section.
\begin{lem}
The set  $J_0 = J^r(W^G) \cap Ker \ \psi$ is a wide subgroupoid of
the groupoid $J^r({\cal D}({\cal C}))$.
\end{lem}
\begin{pf}
Let $a\in G$. Recall that $\Delta(c)=I$ is the constant linear
section. Then $[c]_a$ is the identity at $a$ for $ J^r({\cal
D}({\cal C}))$ and $[c]_a\in J_0$. So $J_0$ is wide in $J^r({\cal
D}({\cal C}))$.

Let $[s ]_a, [t ]_a\in J_0(a, a)$, where $s$ and $t $ are local
linear smooth coadmissible sections with $a\in \Do(s_1)\cap
\Do(t_1)$ and $\alpha (a), \beta (a)\in \Do(s_0)\cap \Do(t_0)$.

Since $J_0 = J^r(W^G) \cap Ker \ \psi$, then we have that
\begin{enumerate}[(i)]
\item $[s ]_a, [t ]_a\in J^r(W^G)$ and so we may assume that the
images of $s$ and $t$ are both contained  in $W^G$ and $s, t$ are
smooth by definition of germs of local linear smooth coadmissible
sections.
\item  $[s ]_a, [t ]_a\in Ker \ \psi$ and this implies that
$\psi ([s ]_a) = \psi ([t ]_a) = 1_a \in {\cal D}({\cal C})$ which
gives $s(a) = t (a) = 1_a$ by definition of the final map.
\end{enumerate}
Therefore $(s(a), t(a))\in W^G\sqcap_{\beta_1} W^G$ and ${\sf
d}(s(a), t(a)) = s(a) -_1 t (a) = 1_a\in W^G$ which implies that
\[   (s(a), t(a))\in (W^G\sqcap_{\beta_1} W^G)\cap {\sf d}^{-1}(W^G).  \]
Since $s$ and $t$ are smooth, then the induced maps
\[    (s_0, t_0) :\Do(s_0)\cap \Do(t_0)\to G\sqcap_{\beta} G
\mbox{ and} \  (s, t) : \Do(s_1)\cap \Do(t_1)\to
W^G\sqcap_{\beta_1} W^G  \] are smooth. But, by condition S3 of
definition \ref{llg} , $(W^G\sqcap_{\beta_1} W^G)\cap {\sf
d}^{-1}(W^G)$ is open in $W^G\sqcap_{\beta_1} W^G$. Since $(G, X)$
is a (globally) Lie groupoid, there exist  open neighbourhoods
$U_1$ of $a$ in $G$, $ U_0$ of $\alpha (a)$, $\beta (a)$ in $ X$
and $\alpha (U_1), \beta (U_1)\subseteq U_0$ such that
\begin{equation*}
  (s, t)(U_1) \subseteq (W^G\sqcap_{\beta_1} W^G)\cap {\sf
d}^{-1}(W^G)  \mbox{   and   }    (s_0, t_0)(U_0)\subseteq
(G\sqcap_{\beta} G)\cap
\partial^{-1}(G),
\end{equation*}
 which implies that $(s, t)(U_1)\subseteq
{\sf d}^{-1}(W^G)$ and $(s_0, t_0)(U_0)\subseteq
{\partial}^{-1}(G)$. Thus $(s * t^{-1})(U_1)\subseteq W^G$ and $(s
* t^{-1})_0(U_0)\subseteq G$, and hence $[s * t^{-1}]_a\in
J^r(W^G)$. Since $s(a) = t(a)$, then $[s * t^{-1} ]_a\in Ker \
\psi.$ Therefore $[s * t^{-1} ]_a\in J_0(a, a)$ and this completes
the proof.
\end{pf}
\begin{lem}\label{nsg}
The groupoid $J_0$ is a normal subgroupoid of the groupoid
$J^r({\cal D}({\cal C}))$.
\end{lem}
\begin{pf}
Let $[k]_a\in J_0 (a, a)$ and let $[s ]_a\in J_0(b, a)$ for some
$a,b \in G$ where $k, s $ are local smooth coadmissible sections,
$\Delta(s)=f$,  with $b = f_1(a) $ and $\beta_1 k (a) = \alpha_1 k
(a) = \beta_1 s (a) = a.$ Moreover $k(a) = 1_a$.  Since $J^r({\cal
D}({\cal C}))$ is generated by $J^r(W^G)$, then
\[   [s ]_a = [s^n]_{a_n} *  \cdots * [s^1]_{a_1}, \ \ \  {s^i}\in \Gamma^r(W^G)  \]
where $a_1 = a$, $a_{i+1} = f_i(a_i)$,   $i = 1, \cdots ,n$.
$[s^i]_{a_i}\in J^r(W^G)$, where we may assume that the images of
the $s^i, i =1,\ldots,n$ are contained in $W^G$ and are smooth.
\begin{align*}
[s]_a [k]_a {[s ]_a}^{-1}  & =  [s^n]_{a_n} *  \cdots *
[s^1]_{a_1}* [k]_a *([s^n]_{a_n} *  \cdots * [s^1]_{a_1})^{-1}
\\
                                    & =  [s^n]_{a_n}* \cdots * [s^1]_{a_1}* [k]_a* {[s^1]_{a_1}}^{-1}  *  \cdots * {[s^n]_{a_n}}^{-1}   \\
                                    & =  [s^n]_{a_n}* \cdots * [s^1]_{a_1}* [k]_a* {[(s^1)^{-1}]_{f_1a_1}} *  \cdots * [(s^n)^{-1}]_{f_1a_n = b} \\
                                    & =  [s * k * s^{-1} ]_b\in J_0(b, b).
\end{align*}
In fact, now, since ${k}^{-1}(a) = -_1k(I^{-1}(a)) = -_1k(a)$,
then ${k}^{-1} (a) = -_1k(a)$. But $k(a) = 1_a$, by definition of
$J_0$; hence $k^{-1}(a) = 1_a\in -_1{W^G}$.

Since, by condition $S1$ of definition \ref{llg}, $W^G =
-_1{W^G}$, then $k(a) \in W^G$. Since $[s ]_a\in J_0(b, a)$, then
we may assume that the image of $s$ is contained in $W^G$ and $s$
is a local linear smooth coadmissible section. So $s(a)\in W^G$,
and therefore
\[ (s(a), -_1k(a))\in W^G\sqcap_{\beta_1} W^G, \ \ \  (s_0(x), -k_0(x))\in G\sqcap_{\beta} G \]
and ${\sf d}(s (a), -k(a)) = s (a) +_1 k(a) = (s * k)(a) = s(a)$.
Also $\partial (s_0(x), -k_0(x)) = s_0(x)+ k_0 (x) = (s_0 *
k_0)(x) = s_0(x)$. Hence $ (s (a), -_1k(a))\in {W^G}_{\sf d}$ and
$ (s_0(x), -k_0(x))\in {G}_{\partial}$, for $x\in X$. By the
smoothness of $k^{-1}$ and $s$, induced map
\[   (s, {k}^{-1} ):\Do({s_1})\cap \Do({k_1}^{-1}) \to  W^G\sqcap_{\beta_1} W^G\mbox{ and } \ (s_0, {k_0}^{-1}
):\Do({s_0})\cap \Do({k_0}^{-1}) \to  G\sqcap_{\beta} G  \] are
smooth. Hence there exists a pair of open neighbourhood $(U_0,
U_1)$ where $\alpha (a), \beta (a) \in U_0$ in $X$, and $a\in U_1$
in $G$ such that
\[ (s, {k}^{-1}) (U_1)\subseteq {W^G}_{\sf d }, \ \  \ (s_0, {k_0}^{-1}) (U_0)\subseteq
{G}_{\partial}  \]
\[ ({s}(U_1)-_1{k}(U_1)) \subseteq W^G, \ \  \ ({s_0}(U_0)-{k_0} (U_0))\subseteq {G}. \]
Therefore $[s * k ]_a\in J_0 (W^G)$.

Thus we may assume that the image of $s * k $ is contained in
$W^G$ and $s * k $ is a local linear  smooth coadmissible section.
Since  $\beta_1 ( s * k) (a) = \beta_1 s (a) = \beta_1 k (a) = a $
and  $(s * k)(a) = s (a)$. Then $((s * k )(a), s (a))\in
W^G\sqcap_{\beta_1} W^G$ and  so ${\sf d} ((s * k )(a), s (a))=
(k* s )(a)-_1 s (a) =1_a\in W^G$, and $((s * k )(a), s (a))\in
{W_{\sf d}}^G$. Similarly  $a\in G(x, y)$, for $x, y\in X$, $((s *
k )_0(x), s_0 (x))\in {G}_{\partial}$. Since $s $ and $ s * k $
are smooth, then they induce  smooth maps
\[   ((s * k ), s): \Do(k_1)\cap \Do(s_1)\to W^G\sqcap_{\beta_1}
W^G, \qquad    ((s * k )_0, s_0): \Do(k_0)\cap \Do(s_0)\to
G\sqcap_{\beta} G. \] But ${W_{\sf d}}^G$ and ${G}_{\partial}$ are
open in $W^G\sqcap_{\beta_1}W^G$ and $G\sqcap_{\beta }G$,
respectively. Hence there exists a pair of neighbourhoods $(U'_0,
U'_1)$ of $\alpha (a), \beta (a)\in U'_0$ in $X$ and $a\in U'_1$
in $G$ such that
\[ ((s * k ), s)(U'_1)\subseteq W^G\sqcap_{\beta_1} W^G, \ \ \
((s * k )_0, s_0)(U'_0)\subseteq G\sqcap_{\beta} G \] which
implies that
\[     (s * k )(U'_1)-_1s (U'_1) \subseteq W^G \  \mbox{\  and  } \ (s * k )_0(U'_0)-s_0 (U'_0) \subseteq G.   \]
Since  $f_1(a) = b$, $[s * k ]_a*{[s ]_a}^{-1} = [s * k]_a
*[{s}^{-1}]_{ b} = [s * k * s^{-1} ]_b\in J(b,b)$. But $[s
* k * s^{-1} ]_b \in (Ker \phi)(b, b)$, since $ (s * k *s^{-1})
(b) = 1_b.$ Hence $[s
* k * s^{-1} ]_b \in J_0(b, b)$ and so $J_0$ is a normal
subgroupoid of $J^r({\cal D}({\cal C}))$.
  \QED  \end{pf}


\section{The Holonomy groupoid of  ${\bf (\D(C), W^G)}$}

In this section, we deal with some locally Lie structures on an
edge symmetric double groupoid  ${\cal D}({\cal C})$ corresponding
to a crossed module ${\cal C}= (C, G, \delta)$ -- namely such a
locally Lie structure is a Lie groupoid structure on the groupoid
$(G, X)$ of ${\cal D}({\cal C})$, and a manifold structure on a
certain subset $W^G$ of the set of squares, satisfying certain
conditions. The Lie groupoid $Hol({\cal D}({\cal C}), W^G)$ we
construct will be called the holonomy groupoid of the $V$-locally
Lie double groupoid $({\cal D}({\cal C}), W^G)$. Further, we state
a universal property of the Lie groupoid $Hol({\cal D}({\cal C}),
W^G)$ in Theorem \ref{gt}.

We state a part of a Lie version of the Brown-Spencer Theorem
given in Brown-Mackenzie \cite{Br-Mac}. We give a definition of
Lie crossed modules of groupoids and of double Lie groupoids.

\begin{defn}{\rm
A {\bf Lie  crossed module of groupoids} is a crossed module $(C,
G, \delta)$ together with a Lie groupoid structure on  $C, G$ (so
that $G,X$ is also a Lie groupoid) such that $\delta :C\to G$ and
the action of $G$ on $C$ are smooth.}
\end{defn}

Recall that a {\bf double groupoid} consists of a quadruple of
sets $(D, H, V, X),$ together with groupoid structures on $H$ and
$V$, both with base $X$, and two groupoid structure on $D$, a
horizontal with base $V$,  and a vertical structure with base $H$,
such that the structure maps (source, target, difference map, and
identity maps) of each structure on $D$ are morphisms with respect
to the other.

\begin{defn}{\rm
A {\bf  double Lie groupoid} is a double groupoid $\cal D$ = $(D;
H, V, X)$ together with differentiable structures on $D$, $H$, $V$
and $X$, such that all four groupoid structures are Lie groupoids
and such that the double source map $D\rightarrow
H\times_{\alpha}V = \{(h, v) : \alpha_1 (h) = \alpha_2 (v) \}$,
$d\mapsto ({\alpha_2}(d), {\alpha_1}(d))$ is  a surjective
submersion, where ${\alpha_2}, {\alpha_1}$ are source  maps on
$D$ vertically and horizontally, respectively.}
\end{defn}

In differential geometry, double Lie groupoids, but usually with
one of the structure totally intransitive, have been considered in
passing by Pradines \cite{Pr2,Pr3}. In general, double Lie
groupoids were investigated by K.Mackenzie in \cite{Mac2} and
Brown and Mackenzie \cite{Br-Mac}.

\begin{thm}{\rm \cite{Br-Mac}} \label{br-sp}
Let ${\cal C} = (C, G, \delta )$ be a Lie crossed module with base
space $X$ and let the anchor map $[ , ]: G\to X\times X$ be
transversal as a smooth map. Then the corresponding edge symmetric
double groupoid  ${\cal D}({\cal C})$ is a  double Lie groupoid.
\end{thm}

We define the quotient groupoid
\[      Hol({\cal D}({\cal C}), W^G) = J^r({\cal D}({\cal C}))/J_0  \]
and call this the {\bf holonomy groupoid} of the locally Lie
groupoid
 $({\cal D}({\cal C}), W^G)$ on ${\cal D}({\cal C})$.

We now state our main theorem.

\begin{thm}\label{gt}
Let ${\cal C} = (C, G, \delta)$  be a crossed module and let
${\cal D}({\cal C})$ be the corresponding double groupoid. Let
$({\cal D}({\cal C}), W^G)$ be a $V$-locally Lie double groupoid
for the double groupoid ${\cal D}({\cal C})$. Then there is a Lie
groupoid structure on $Hol({\cal D}({\cal C}), W^G)$, a morphism
\[             \psi : Hol({\cal D}({\cal C}), W^G) \to {\cal D}({\cal C})  \]
of groupoids, and an embedding
\[    i: W^G\to Hol({\cal D}({\cal C}), W^G)  \]
of $W^G$ to an open neighbourhood of $Ob(Hol({\cal D}({\cal C}),
W^G)) = G$, such that
\begin{enumerate}[\rm i)]
\item $\psi $ is the identity on $G$, $\psi i$ is the identity on
$W^G$, ${\psi }^{-1}(W^G)$ is open in  $Hol({\cal D}({\cal C}),
W^G)$, and the restriction $\psi_{W^G} : {\psi }^{-1}(W^G)\to W^G
$ of $\psi$  is smooth.

\item if ${\cal A} = (A, B, \delta')$ is a Lie crossed module and
$\mu :{\cal D}({\cal A}) \to {\cal D}({\cal  C})$ is a morphism of
groupoids such that
\end{enumerate}
\begin{enumerate}[ a)]
\item $\mu$ is the identity on objects;

\item the restriction $\mu_{W^G}:{\mu }^{-1}(W^G)\to W^G$ is smooth and
${\mu }^{-1}(W^G)$ is open in ${\cal D}({\cal A})$ and generates
$\Do({\cal A})$ with respect to $+_1$ as a groupoid.

\item the triple $(\alpha_1, \beta_1, {\cal D}({\cal A}))$
has enough local linear smooth coadmissible sections;
\end{enumerate}
then there is a unique morphism $\mu': ({\cal D}({\cal A}), B) \to
{Hol({\cal D}({\cal C}), W^G)}$ of Lie groupoids such that $\psi
\mu' = \mu$ and $ \mu' (w) = (i \mu )(a)$ for $w\in {\mu
}^{-1}(W^G)$.
\end{thm}
The proof occupies the next two sections.


\section{Lie groupoid structure on $Hol({\cal D}({\cal C}), W^G)$}

The aim of this section is to construct a topology  on the
holonomy groupoid $Hol({\cal D}({\cal C}), W^G)$ such that
$Hol({\cal D}({\cal C}), W^G)$ with this topology is a Lie
groupoid. In the next section we verify that the universal
property of   Theorem \ref{gt} holds. The intuition is that first
of all $W^G$ embeds in $Hol({\cal D}({\cal C}), W^G)$, and second
that $Hol({\cal D}({\cal C}), W^G)$ has enough local linear
coadmissible sections for it to obtain a topology by translation
of the topology of $W^G$.

Let $s\in \Gamma ^r({\cal D}({\cal C}), W^G)$. We define a partial
function $\chi_{s}:W^G\to Hol({\cal D}({\cal C}), W^G)$. The
domain of $\chi_{s}$ is the set of $w\in W^G$ such that $\alpha_1
(w) = a\in \Do(s_1)$ and $\alpha (a), \beta (a)\in \Do(s_0)$. The
value $\chi_{s}(w)$ is obtained as follows. Choose a local linear
smooth coadmissible section $\theta$ through $w$. Then we set
\[    \chi_{s}(w) =\<s \>_{\alpha_1 (w)} \<\theta\>_{\beta_1 (w)}
= \<s * \theta\>_{\beta_1 (w)}. \]

We have to show that $\chi_{s }(w)$ is independent of the choice
of the local linear smooth coadmissible section $\theta $. For
this reason we state a lemma.
\begin{lem}\label{51}
Let $w\in W^G$, and let $s$ and $t$ be local linear smooth
coadmissible sections through $w$. Let $a = \beta_1 w$. Then
$\<s\>_a = \<t\>_a$ in $Hol({\cal D}({\cal C}), W^G)$.
\end{lem}
\begin{pf}
By assumption $s a = t a = w$. Let $b = \alpha_1 w$. Without loss
of generality, we may assume that $s$ and $t$ have the same domain
$(U_0, U_1)$ and have image contained in $W^G$ and $G$,
respectively. By Lemma \ref{st}, $s * t^{-1} \in \Gamma^r (W^G)$.
So $[s * t^{-1}]_b\in J_0$. Hence
\[     \<t \>_a = \<s * t^{-1}\>_b\<t\>_a   = \<s * t^{-1} *t\>_a =  \<s \>_a   \]
\end{pf}

\begin{lem}
$\chi_{s}$ is injective.
\end{lem}
\begin{pf}
Suppose  $\chi_{s} v = \chi_{s} w$. Then $\beta_1 w = \beta_1 v =
a$, say and $ \alpha_1 s \alpha_1 v = \alpha_1 s \alpha_1 w.$ By
definition of $s$, $\alpha_1 v = \alpha_1 w = d$, say.  Let
$\theta$ be a local linear smooth coadmissible section through
$w$. Then we now obtain from $\chi_{s}v =  \chi_{s} w$ that
\[  \<s \>_d    \< \theta \>_a  =\<s \>_d \<\theta'\>_a     \]
and hence, since $Hol({\cal D}({\cal C}), W^G)$ is a groupoid,
that $\<\theta \>_a = \<\theta'\>_a$. Hence $ v = \theta (a) =
{\theta }' (a) = w\in W^G$.
\end{pf}

Let $s\in \Gamma ({\cal D}({\cal C})) $. Then $s$ defines a left
translation $L_{s }$ on ${\cal D}({\cal C})$ by
\[       L_{s} (w) = s (\alpha_1 (w)) +_1 w.  \]
This is an injective partial function on ${\cal D}({\cal C})$. The
inverse ${L_s}^{-1}$ of $L_{s}$ is
\[  {v\mapsto -_1 s (\alpha_1 s)^{-1}(\alpha_1 (v))+_1 v} \]
and ${L_{s }}^{-1} = L_{s^{-1}}$. We call $L_{s }$ the {\bf left
translation corresponding to $s$}.

So we have an injective function $\chi_{s}$ from an open subset of
$W^G$ to $Hol({\cal D}({\cal C}), W^G)$. By definition of
$Hol({\cal D}({\cal C}), W^G)$, every element of $Hol({\cal
D}({\cal C}), W^G))$ is in the image of $\chi_{s}$ for some $s$.
These $\chi_{s}$ will form a set of charts and so induce a
topology on $Hol({\cal D}({\cal C}), W^G)$. The compatibility of
these charts results from the following lemma, which is essential
to ensure that $W^G$ retains its topology in $Hol({\cal D}({\cal
C}), W^G)$ and is open in $Hol({\cal D}({\cal C}), W^G)$. As in
the groupoid case \cite{Aof-Br}, this is a key lemma.

\begin{lem}
Let $s, t \in \Gamma^r ({\cal D}({\cal C}),  W^G)$. Then $
(\chi_{t})^{-1} \chi_{s}$ coincides with $L_{\eta},$ left
translation by the local linear smooth coadmissible section $\eta
=t^{-1} * s $, and $L_{\eta}$ maps open sets of $W^G$
diffeomorphically to open sets of $W^G$.
\end{lem}
\begin{pf}
Suppose $v, w\in W^G$ and $\chi_{s} v = \chi_{t} w$. Choose  local
linear smooth coadmissible  sections $\theta$ and $\theta'$
through $v$ and $w$ respectively such that the images of  $\theta$
and $\theta'$ are contained in $W^G$. Since  $\chi_{s} v =
\chi_{t} w$, then $\beta_1 v = \beta_1 w = a$ say. Let $\alpha_1 v
= b$, $\alpha_1 w = c$.

Since $\chi_{s} v = \chi_{t} w$, we have
\[    \<s * \theta \>_a = \< t * \theta' \>_a     \]
Hence there exists a local linear smooth coadmissible section
$\zeta $ with $a\in \Do(\zeta)$ such that $[\zeta ]_a\in J_0$ and
\[     [s * \theta ]_a =  [t*\theta']_a [\zeta ]_a    \]
Let $\eta ={t }^{-1}* s $. Then in the semigroup $\Gamma^r({\cal
D}({\cal C}), W)$ we have from the above that $\eta *\theta =
\theta' *\zeta$ locally near $a$. So we get $ w = (\theta'*\zeta)(
a) =  {\theta'}(a)+_1 \zeta (a) =
 {\theta'}(a)+_1 1_a = (\eta * \theta) a = \eta \alpha_1 v +_1  v$. This
shows that $(\chi_{t})^{-1}\chi_{s} = L_{\eta}$, left translation
by $\eta \in \Gamma ({\cal D}({\cal C}))$, i.e.,
\begin{align*}
(\chi_{t})^{-1} (\chi_{s }) (v) & =  (\chi_{t})^{-1}(\<s * \theta
\>_{\beta_1 v = a}) \\
                              & =  (t^{-1} * s * \theta )(a), \\
                              & =  (\eta * \theta )(a) , \  \   \mbox{since} \
                               \  \eta = t^{-1}*s  \\
                              & =  \eta (\alpha_1 (\theta (a)) +_1\theta (a), \ \  \mbox{by definition of} \ \ * \\
                              & =   \eta (\alpha_1 (v))+_1 v, \ \  \mbox{since}  \ \ \ \theta(a) = v \\
                              & =  L_{\eta } (v), \ \ \  \mbox{by definition of} \ \ \ L_{\eta}.
\end{align*}
However, we also have $\eta = \theta' *  \zeta * {\theta}^{-1}$
near $\alpha_1 v$. Hence $L_{\eta } =L_{\theta'} L_{\zeta }
L_{{\theta}^{-1}}$ near $v$. Now $L_{{\theta}^{-1}}$ maps $v$ to
$1_a$, $L_{\zeta}$ maps $1_a$ to $1_a$, and $L_{\theta'}$ maps
$1_a$ to $w$. We prove the first of these the others being
similar.
\begin{align*}
L_{{\theta}^{-1}}(v) & =   {\theta}^{-1}(\alpha_1 (v)) +_1 v \\
                     & =  -_1\theta(\alpha_1 \theta)^{-1}(\alpha_1 v) +_1 v, \ \  \mbox{by definition of} \ \ {\theta}^{-1} \\
                     & =  -_1\theta(\beta_1(v))+_1 \theta(\beta_1 v), \ \ \mbox{since}\ \
                     \theta(\beta_1 v) = v \\
                     & =  1_a,
\end{align*}
So these left translations are defined and smooth on open
neighbourhoods of $v$, $1_a$ and $1_a$ respectively. Hence
$L_{\eta }$ is defined and smooth on an open neighbourhood of $v$.
\end{pf}

We now impose on $Hol({\cal D}({\cal C}), W^G)$ the initial
topology with respect to the charts $\chi_{s}$ for all $s \in
\Gamma^r ({\cal D}({\cal C}), W^G)$. In this topology each element
$h$ has an open neighbourhood diffeomorphic to an  open
neighbourhood of $1_{\beta_1 h}$ in $W^G$.
\begin{lem}
With the above topology, $Hol({\cal D}({\cal C}), W^G)$ is a Lie
groupoid.
\end{lem}
\begin{pf}
Source and target maps are smooth: In fact, for $w\in W^G$,
\[  \beta_H (\chi_{s } (w)) = \beta_1 (w), \ \ \
\alpha_H (\chi_{s }(w)) = \alpha_1 (s \alpha_1 (w)).  \] It
follows that $\alpha_H$ and $\beta_H$ are smooth.

Now we have to prove that
\[     {\sf d}_H :Hol({\cal D}({\cal C}), W^G)\sqcap_{\beta } Hol({\cal D}({\cal C}), W^G)\to Hol({\cal D}({\cal C}), W^G) \]
is smooth. Let $\<s \>_a, \<t \>_a\in Hol({\cal D}({\cal C}),
W^G)$. Then $\chi_{s }(1_a) = \<s \>_a, \chi_{t} (1_a) = \<t
\>_a$, and if $\eta = {s }* t^{-1}$, then $\chi_{\eta}(1_b) = \<{s
}*t^{-1} \>_b$ where $b = \beta_1 t (a)$. Let $v\in \Do(\chi_{s}),
w\in \Do(\chi_{t})$, with $\beta_1 v = \beta_1 w = c$, say and let
$\theta $ and $\theta'$ be elements of $\Gamma^r(W^G)$ through $v$
and $w$ respectively. Let $d = \beta_1 ( t * \theta' ) (c).$ Then
\begin{align*}
{\chi_{\eta}}^{-1} {\sf d}_H (\chi_{s} \times \chi_{t })(v, w) & =
{\chi_{\eta}}^{-1} {\sf d}_H (\chi_{s }(v), \chi_{t}(w)) \\
              & =  {\chi_{\eta}}^{-1} {\sf d}_H(\<s * \theta\>_c, \<t * \theta' \>_c),\ \ \mbox{by definition of} \ \chi_s, \chi_t \\
              & =  \chi_{\eta}^{-1}(\<(s * \theta ) * (t * \theta')^{-1} \>_d), \ \ \mbox{by definition of} \  {\sf d}_H \\
              & =  ( {\eta}^{-1})*(s * \theta ) * (t * \theta' )^{-1}(d), \ \ \mbox{by definition of} \  {\chi_{\eta}}^{-1} \\
              & =  ({({s}* {t }^{-1})}^{-1})*(s * \theta ) * (t * \theta' )^{-1}(d), \ \ \mbox{since}\  \eta = (s*t^{-1}) \\
              & =  ( t*s^{-1}* s * \theta ) * (t * \theta' )^{-1})(d) \\
              & =  ((t * \theta ) * (t * \theta')^{-1})(d) \\
              & =  ((t * \theta )(\alpha_1  (t * \theta')^{-1}(d) +_1 (t*\theta')^{-1}(d) \\
              & =  (t * \theta )_1(c )-_1 (t *\theta' )(\alpha_1 (t*\theta')^{-1})(d), \ \ \mbox{since}\
              \alpha_1(t*\theta')^{-1}(d)=c \\
              & =   t (\alpha_1 \theta_1 (c)+_1 \theta (c)-_1(t
              (\alpha_1 \theta'(c) +_1 \theta'(c))   \\
              & =  (t (\alpha_1 (v))+_1 v -_1 (t (\alpha_1 (w)) +_1 w) \\
              & =   L_{t }(v) -_1 L_{t }(w)  \\
              & =  \Omega (v, w),
\end{align*}
say. The smoothness of this map $\Omega$ at $(1_a, 1_a)$ is now
easily shown by writing $t = t_n * \cdots* t_1$ where $t_i\in
\Gamma^r(W^G)$ and using induction and a similar argument to  that
of Lemma \ref{nsg}.
\end{pf}


\section{The Universal Property of Hol({\cal D}({\cal C}), $W^G)$ }

In this section we state and prove the main theorem of the
universal property of the morphism $\psi: Hol({\cal D}({\cal C}),
W^G) \to {\cal D}({\cal C}).$ Note that for the case of groupoids
rather than crossed modules, Pradines \cite{Pr1} stated a
differential version involving germs of locally Lie groupoids in
\cite{Pr1}, and formulated the theorem in terms of adjoint
functors. No information was given on the construction or proof. A
version for locally topological groupoids was given in Aof-Brown
\cite{Aof-Br}, with complete details of the construction and
proof, based on conversations of Brown with Pradines. The
modifications for the differential case were given in Brown-Mucuk
\cite{Br-Muc2}.

The main idea is when we are given a $V$-locally Lie double
groupoid $({\cal D}({\cal C}), W^G)$ of a double groupoid ${\cal
D}({\cal C})$, coming from a Lie crossed module ${\cal C}$, a Lie
crossed module ${\cal A}$ and a morphism
\[ \mu :{\cal D}({\cal A})\to {\cal D}({\cal C}) \]
with suitable conditions, we can construct a morphism
\[    \mu' : {\cal D}({\cal A})\to Hol({\cal D}({\cal C}), W^G), \]
where $Hol({\cal D}({\cal C}), W^G))$ is a holonomy groupoid of a
locally Lie crossed module, such that
\[    \phi \mu' = \mu.  \]
We prove that such a morphism $\mu'$ is well-defined, smooth and
unique. Now let $({\cal D}({\cal C}), W^G)$ be a $V$- locally Lie
double groupoid as above.

\begin{thm} \label{hol-univ2}
If ${\cal A} = (A, B, \delta')$ is a Lie crossed module and $\mu
:{\cal D} ({\cal A})\to {\cal D}({\cal C})$ is a morphism of
groupoid such that
\begin{enumerate}[\rm (i)]
\item $\mu$ is the identity on objects;

\item the restriction $\mu_{W^G} : {\mu }^{-1}(W^G)\to W^G$ of $\mu$ is smooth
and $\mu^{-1}(W^G)$ is open in ${\cal D}({\cal A})$ and generates
${\cal D}({\cal A})$ as a groupoid with respect to $+_1$.

\item  the triple $(\alpha_1, \beta_1, {\cal D}({\cal A}))$ has enough local smooth
coadmissible sections.
\end{enumerate}
Then there exists a unique morphism
\[    \mu' : {\cal D}({\cal A})\to Hol({\cal D}({\cal C}), W^G) \]
of Lie groupoids such that $\psi \mu' = \mu$ and $\mu' (w) = i \mu
(w)$ for $w\in \mu^{-1}(W^G)$.
\end{thm}
\begin{pf}
Since, by condition {\rm (i)}, $\mu$ is the identity on $G$, then
$G = B$ and $X = X'$ which implies that $\mu (G) = G$, $\mu (X) =
X$. But $G\subseteq W^G\subseteq {\cal D}({\cal C})$, by condition
$S2$ of Definition \ref{llg}. Hence $\mu (G) \subseteq
W^G\subseteq {\cal D}({\cal C})$. So it follows that $G\subseteq
\mu^{-1}(W^G)\subseteq {\cal D}({\cal A})$.

But, by condition {\rm (ii)}, $\mu^{-1}(W^G)$ is an open in ${\cal
D}({\cal A})$. Hence $\mu^{-1}(W^G)$ is open neighbourhood of $G$
in ${\cal D}({\cal A})$. Since $\mu^{-1}(W^G)$ generates
$\Do({\cal A})$, we can write $w = w_n +_1\cdots +_1 w_1$, where
$\mu (w_i)\in W^G$,  $i =1,\ldots,n$.

Since $(\alpha_1, \beta_1, {\cal D}({\cal A}))$ has enough local
linear smooth coadmissible sections, by condition {\rm (iii)}, we
can choose local linear smooth coadmissible sections $\theta_i$
through $w_i$, $i = 1,\ldots, n$, such that they are composable
and their images are contained in $\mu^{-1}(W^G)$.

Because of the condition {\rm (ii)}, the smoothness of $\mu$ on
$\mu^{-1}(W^G)$ implies that $\mu \theta_i$ is a local linear
smooth coadmissible section through $\mu (w_i)\in W^G$ whose image
is contained in $W^G$. Therefore $\mu \theta \in \Gamma^r({\cal
D}({\cal C}), W^G)$. Hence we can set
\[         \mu' (w) = \<\mu \theta \>_{\beta_1 (w)}   \]
\end{pf}
\begin{lem}
$\mu' (w)$ is independent of the choices which have been made.
\end{lem}
\begin{pf}
Let $ w = v_m +_1 \cdots +_1 v_1$, where  $\mu v_j\in W^G$ and $ j
= 1,\cdots,m$. Choose a set of local linear smooth coadmissible
sections ${\theta'}_j$ through $v_j$ such that the ${\theta'}_j$
are composable and their images are contained in $\mu^{-1}(W^G)$.

Let $\theta' = {\theta'}_m * \cdots * {\theta'}_1$. Then $\mu
\theta'\in \Gamma^r({\cal D}({\cal C}), W^G)$, and so $\<\mu
\theta'\>_c\in Hol({\cal D}({\cal C}), W^G)$. Since by assumption,
$\theta (c) = \theta' (c) = w\in {\cal D}({\cal A})$, then
$(\theta (c), \theta'(c))\in {\cal D}({\cal A})\sqcap_{\beta_1 }
{\cal D}({\cal A})$ and ${\sf d}_A (\theta (c), \theta' (c)) =
{\theta}(c) -_1 \theta'(c) = 1_c.$ Hence $(\theta (c),
\theta'(c))\in {{\sf d}_A}^{-1} \mu^{-1}(W^G)$ because $1_c \in
\mu^{-1}(W^G)$.

Because ${\cal A}$ is a Lie crossed module and the corresponding
double groupoid ${\cal D}({\cal A})$ is a  double Lie groupoid,
the difference map ${\sf d}_A : {\cal D}({\cal A})\sqcap_{\beta_1
} {\cal D}({\cal A})\to {\cal D}({\cal A})$ is smooth. Since
$\mu^{-1}(W^G)$ is open in $ {\cal D}({\cal A})$, by condition
{\rm (ii)}, then ${{\sf d}_A}^{-1} \mu^{-1}(W^G)$ is open in
${\cal D}({\cal A})\sqcap_{\beta_1 }{\cal D}({\cal A})$.

But, by the smoothness of $\theta $ and $\theta'$, the induced
maps $(\theta, \theta') : \Do(\theta_1 )\cap \Do(\theta'_1)\to
{\cal D}({\cal A})\sqcap_{\beta_1} {\cal D}({\cal A})$,
$(\theta_0, \theta'_0) : \Do(\theta_0 )\cap \Do(\theta'_0)\to G
\sqcap_\beta G$ are smooth. Hence there exists open neighbourhoods
$N$ of $c$ in $G$ and $N_0$ of both $\alpha (c), \beta (c)$ such
that $(\theta, {\theta}') (N) \subseteq ({{\sf d}_A }^{-1}{\mu
}^{-1})(W^G)$ and $(\theta_0, \theta'_0) (N_0) \subseteq
\partial_A^{-1}(G)$. This implies that ${\theta}* {{\theta'}}^{-1}(\alpha_1
{\theta'} N)\subseteq {\mu }^{-1}(W^G)$ and $\theta_0
*{\theta_0'}^{-1}(N_0) \subseteq G$. So, after suitably
restricting $\theta$, $\theta'$, which we may suppose done without
change of notation, we have that $\theta
* {\theta'}^{-1}$ is a local linear smooth coadmissible section
through $1_d\in {\cal D}({\cal A})$ and its image is contained in
${\mu }^{-1}(W^G)$. So $\mu (\theta
* {\theta'}^{-1})$ is a local linear smooth coadmissible section
through $1_d\in W^G$, and its image is contained in $W^G$.
Therefore  $[ \mu (\theta * {\theta'}^{-1})]_d\in J^r(W^G)$.

Since $\theta (c) = {\theta'}(c)$, then $\psi [\mu \theta ]_c =
\psi [\mu \theta' ]_c$. But $\psi $ and $\mu$ are morphisms of
groupoids; hence $\psi[\mu (\theta * {\theta'}^{-1})]_d = 1_d$,
and so $[\mu  (\theta * {\theta'}^{-1})]_d\in Ker \psi.$ Therefore
$[\mu (\theta * {\theta'}^{-1})]_d\in J^r(W)\cap Ker \psi = J_0$.
Since $\mu $ is a morphism of groupoids, we have $[\mu (\theta *
{\theta'}^{-1})]_d\in J^r.$ Hence  $\< \mu (\theta *
{\theta'}^{-1})\>_d = 1_d \in Hol({\cal D}({\cal C}), W^G))$, and
so
\[ \< \mu \theta\>_c = \< \mu  \theta\>_c \<\mu (\theta * {\theta'}^{-1}) \>_d = \< \mu \theta'\>_c \]
which shows that $\mu' w$ is independent of the choices made.
\end{pf}
\begin{lem}
$\mu'$ is a morphism of groupoids.
\end{lem}
\begin{pf}
Let $ u = w +_1 v$ be an element of ${\cal D}({\cal A})$ such that
$ w = w_n +_1 \cdots +_1 w_1$  and  $v = v_m +_1 \cdots +_1 v_1$,
where $w_i, v_j\in {\mu }^{-1}(W^G)$, $i =1,\ldots,n$ and $j =
1,\ldots,m$. Then $ u = w_n +_1 \cdots +_1 w_1 +_1 v_m +_1 \cdots
+_1 v_1$.

Let $\theta_i, {\theta'}_j$ be local linear smooth coadmissible
section through $w_i$ and $v_j$ respectively such that they are
composable and their images are contained in ${\mu }^{-1}(W^G)$.
Let $\theta = \theta_n * \cdots * \theta_1$ and $\theta' =
{\theta'}_m * \cdots * {\theta'}_1$, $\kappa = \theta * \theta'$.
Then $\kappa$ is a local linear smooth coadmissible section
through $u\in {\cal D}({\cal A})$, and $\mu \theta$, $ \mu
\theta'$, $\mu \kappa\in \Gamma^r({\cal D}({\cal C}), W^G)$, and
$\mu \kappa = \mu \theta * \mu \theta'$, since $\mu$ is  a
morphism of groupoids.

Let $ a = \beta_1 w$, $ b = \beta_1 v$. Then $\< \mu \kappa \>_a =
\<\mu \theta \>_a \<\mu \theta \>_b$ and so $\mu'$ is a  morphism.
\end{pf}
\begin{lem}
The morphism $\mu'$ is smooth, and is the only morphism of
groupoids such that $\psi \mu' = \mu$ and $\mu' a = (i \mu )(a)$
for all $a\in {\mu }^{-1}(W^G)$.
\end{lem}

\begin{pf}

Since $(\alpha_1, \beta_1, {\cal D}({\cal A}))$ has enough  local
linear smooth coadmissible section, it is enough to prove that
$\mu'$ is smooth at $1_a$ for all $a\in G$. Let {\bf c} denote the
linear  coadmissible section ${\bf c} :G\to {\cal D}({\cal C})$,
$a\mapsto 1_a$ and $c_0:X\to G$, $x\mapsto 1_x$.

Let $a\in G$. If $w\in {\mu}^{-1}(W^G)$ and $s $ is a local linear
smooth coadmissible section through $w$, then $\mu' w = \< \mu s
\>_{\beta_1 w} = \chi_{\bf c} \mu (w)$.  Since $\mu$ is smooth, it
follows that $\mu'$ is smooth.

The uniqueness of $\mu'$ follows from the fact that $\mu'$ is
determined on ${\mu }^{-1}(W^G)$ and that this set generates
${\cal D}({\cal A})$.
\end{pf}
This completes the proof of our main result, Theorem \ref{gt}.


{}

\end{document}